\documentclass[12pt,a4paper]{article}

\usepackage{amsmath}
\usepackage{amsfonts}
\usepackage{amssymb}
\usepackage{latexsym}
\usepackage{graphicx}

\usepackage{dsfont}

\makeatletter
\@addtoreset{equation}{section}
\makeatother

\newcommand{\be}{\begin{equation}}
\newcommand{\ee}{\end{equation}}
\newcommand{\nnr}{\nonumber \\}

\newcommand{\eq}[1]{(\ref{#1})}
\newcommand{\fr}{\frac}
\newcommand{\tf}{\tfrac}

\newcommand{\df}{\textrm{d}}
\newcommand{\expe}[1]{\textrm{e}^{#1}}
\newcommand{\pd}{\partial}

\newcommand{\Real}{\textrm{Re} \,}

\newcommand{\bb}[1]{\mathbb{#1}}

\renewcommand{\vec}[1]{\mathbf{#1}}
\newcommand{\gvec}[1]{\boldsymbol{#1}}

\newcommand{\Var}{\textrm{Var}}
\newcommand{\Cov}{\textrm{Cov}}
\newcommand{\Corr}{\textrm{Corr}}

\begin{document}


\thispagestyle{empty}

\begin{center}
\textbf{\Large{Schl\"{o}milch integrals and probability distributions on the simplex}}\\
\vspace{50pt}
\large{David D. K. Chow}
\end{center}

\begin{center}
{\bf Abstract\\}
\end{center}
The Schl\"{o}milch integral, a generalization of the Dirichlet integral on the simplex, and related probability distributions are reviewed.  A distribution that unifies several generalizations of the Dirichlet distribution is presented, with special cases including the scaled Dirichlet distribution and certain Dirichlet mixture distributions.  Moments and log-ratio covariances are found, where tractable.  The normalization of the distribution motivates a definition, in terms of a simplex integral representation, of complete homogeneous symmetric polynomials of fractional degree.

\newpage


\section{Introduction}


One of the most important multivariate integrals is the Dirichlet integral on the simplex \cite{dirichlet}, which is a multivariate generalization of the Euler beta function.  Expressed as an integral over the simplex $S_K = \{ \vec y = (y_1 , \ldots , y_K)  \subset \bb R^K : \sum_{i = 1}^K y_i = 1 , y_i \geq 0 \}$, the Dirichlet integral is
\be
\int_{S_K} \df^{K - 1} y \, \! \prod_{i = 1}^K y_i^{\alpha_i - 1} = \fr{\prod_{i = 1}^K \Gamma (\alpha_i)}{\Gamma (\alpha_+)} ,
\label{dirichletint}
\ee
where $\alpha_i > 0$ and $\alpha_+ = \sum_{i = 1}^K \alpha_i$.  Because of the sum constraint $\sum_{i = 1}^K y_i = 1$, we may regard the integrand as a function of $y_1 , \ldots , y_{K - 1}$ and take
\be
\int_{S_K} \df^{K - 1} y = \int_0^1 \! \df y_1 \, \int_0^{1 - y_1} \! \df y_2 \ldots \int_0^{1 - y_1 - \ldots - y_{K - 2}} \! \df y_{K - 1} .
\ee

A generalization is the Schl\"{o}milch integral
\be
\int_{S_K} \df^{K - 1} x \, \! \fr{\prod_{i = 1}^K x_i^{\alpha_i - 1}}{(\sum_{j = 1}^K \beta_j x_j)^{\alpha_+}} = \fr{\prod_{i = 1}^K \Gamma (\alpha_i)}{\Gamma (\alpha_+) \prod_{i = 1}^K \beta_i^{\alpha_i}} ,
\label{schlomilchint}
\ee
where $\beta_i > 0$, which was obtained by Schl\"{o}milch not long after the discovery of the Dirichlet integral, appearing first in the textbook \cite{schlomilch0} (p.\ 165), and then further developed in the article \cite{schlomilch}.  The Schl\"{o}milch integral \eq{schlomilchint} follows from the Dirichlet integral \eq{dirichletint} through the coordinate transformation
\be
y_i = \fr{\beta_i x_i}{\sum_{j = 1}^K \beta_j x_j} ,
\label{transform}
\ee
although it was originally derived in a more complicated manner.  If all $\beta_i$ are equal, then this transformation is trivial and the Schl\"{o}milch integral \eq{schlomilchint} reduces to the Dirichlet integral \eq{dirichletint}.

Schl\"{o}milch gave a slightly more general result,
\be
\int_{x_i \geq 0, \, \lambda \leq \sum_{i = 1}^K x_i \leq \kappa} \! \df^K x \, \fr{\prod_{i = 1}^K x_i^{\alpha_i - 1}}{(\theta + \sum_{j = 1}^K \gamma_j x_j)^{\alpha_+}} f({\textstyle \sum_{k = 1}^K x_k}) = \fr{\prod_{i = 1}^K \Gamma (\alpha_i)}{\Gamma (\alpha_+)} \int_\lambda^\kappa \! \df t \, \fr{t^{\alpha_+ - 1} f(t)}{\prod_{j = 1}^K (\theta + \gamma_j t)^{\alpha_j}} ,
\label{schlomilchint2}
\ee
where $f$ is an arbitrary function, and we take $\beta_i = \theta + \gamma_i$ to facilitate comparison with the original literature.  This type of generalization is analogous to Liouville's well-known generalization of the Dirichlet integral.  Taking $f(t) = \delta(t - 1 + \epsilon)$, i.e.\ a Dirac delta function centred at $1 - \epsilon$ for some $\epsilon > 0$, as well as $\lambda < 1 < \kappa$, then taking $\epsilon \rightarrow 0$, we obtain the Schl\"{o}milch integral \eq{schlomilchint}.  If all $\gamma_i$ are equal, then the more general Schl\"{o}milch integral \eq{schlomilchint2} reduces to Liouville's generalization of the Dirichlet integral.

The Schl\"{o}milch integral \eq{schlomilchint} or the slightly more general \eq{schlomilchint2} have appeared in widely used textbooks from the late 19th and early 20th centuries, such as Todhunter (pp.\ 263--265 of \cite{todhunter}) and Edwards (pp.\ 169--173 of \cite{edwards}).  The result does not seem to have been widely used in research at that time, but an example is Dixon \cite{dixon}, who does not provide a reference, instead providing a self-contained derivation.

In contemporary use, the Schl\"{o}milch integral has applications in several areas of science and mathematics, and appears in reference books such as Gradshteyn and Ryzhik (entry 4.637 of \cite{graryz}), Prudnikov, Brychkov and Marichev (Section 3.3.4 of \cite{prbrma}), and Zwillinger (p.\ 102 of \cite{zwillinger}).  The formula \eq{schlomilchint} is useful in both directions: it can be used directly as an evaluation of an integral on a simplex, but can also be used, as an intermediate step in a larger calculation, to replace $\prod_{i = 1}^K \beta_i^{-\alpha_i}$ with a simplex integral.  Perhaps most notably, the Schl\"{o}milch integral is in widespread use in quantum field theory, where it is a key tool for evaluating Feynman diagrams, with $x_1 , \ldots , x_{K - 1}$ referred to as Feynman parameters, and the Schl\"{o}milch integral \eq{schlomilchint} sometimes attributed to Feynman.  In this context, Feynman stated that the technique was suggested by Schwinger \cite{feynman}, as further explained in the historical study \cite{schweber} (see in particular pp.\ 445 and 452--454).  For a standard textbook treatment, see e.g.\ equation (6.42) of \cite{pessch}; a recent review of Feynman integrals \cite{weinzierl} considers Feynman parameters in Section 2.5.3.  Despite its widespread use, the Schl\"{o}milch integral is not as well-known as it should be, as suggested by its rediscovery over the last few years in disparate fields \cite{xu, cocmcd}; the Dirichlet integral itself has also been rederived recently \cite{verseg}.  Furthermore, the origin of the Schl\"{o}milch integral appears to have been completely forgotten. 

Schl\"{o}milch pointed out that further results can be obtained by differentiating with respect to parameters.  In equation (6) of \cite{schlomilch}, by differentiating with respect to $\theta$, is an equation equivalent to
\be
\int_{\substack{\\
\\
x_i \geq 0, \\
0 \leq \sum_{i = 1}^K x_i \leq 1}}
\! \df^K x \, \fr{\prod_{i = 1}^K x_i^{\alpha_i - 1} f(\sum_{k = 1}^K x_k)}{(\theta + \sum_{j = 1}^K \gamma_j x_j)^{\alpha_+ + 1}} = \fr{\prod_{i = 1}^K \Gamma (\alpha_i)}{\Gamma (\alpha_+ + 1)} \int_0^1 \! \df t \, \fr{t^{\alpha_+ - 1} f(t)}{\prod_{j = 1}^K (\theta + \gamma_j t)^{\alpha_j}} \sum_{i = 1}^K \fr{\alpha_i}{\theta + \gamma_i t} .
\label{schlohigher}
\ee
Again by taking $f(t) = \delta(t - 1 + \epsilon)$ and then the limit $\epsilon \rightarrow 0$, we obtain the simplex integral
\be
\int_{S_K} \df^{K - 1} x \, \! \fr{\prod_{i = 1}^K x_i^{\alpha_i - 1}}{(\sum_{j = 1}^K \beta_j x_j)^{\alpha_+ + 1}} = \fr{\prod_{i = 1}^K \Gamma (\alpha_i)}{\Gamma (\alpha_+ + 1) \prod_{j = 1}^K \beta_j^{\alpha_j}} \sum_{i = 1}^K \fr{\alpha_i}{\beta_i} .
\label{schloint1}
\ee
This result may also be obtained by differentiating the Schl\"{o}milch integral \eq{schlomilchint} with respect to each $\beta_i$ and then summing.  The differentiation procedure increases the exponent of the denominator of the integrand from $\alpha_+$ to $\alpha_+ + 1$ and can be repeated, increasing the exponent to $\alpha_+ + n$, for positive integers $n$.  In this article, we shall find explicit expressions for such integrals, and examine some further generalizations.

The Dirichlet integral and its generalizations have wide-ranging applications, but we shall focus here on applications to associated probability distributions on the simplex, for which these integrals appear in normalization constants.  Through a special mixture distribution based on the Schl\"{o}milch integral, we shall provide a unified framework for several distributions appearing in the literature.  Special cases include the G3D generalized Dirichlet distribution \cite{chenov} or scaled Dirichlet distribution \cite{aitchison}, the shifted-scaled Dirichlet distribution \cite{momapaeg} or simplicial generalized beta distribution \cite{graf}, the flexible Dirichlet distribution \cite{onmimo}, and the tilted Dirichlet distribution \cite{coltaw}.  In a similar manner, we shall also obtain a family of distributions that includes the Concrete or Gumbel-softmax distribution \cite{mamnte, jagupo}.  By including additional parameters, these generalizations of the Dirichlet distribution allow for behaviours that the Dirichlet distribution itself cannot model, for example positive covariances $\Cov(X_i, X_j)$.

Although we can show, through explicit computation, examples in which $\Cov (X_i, X_j)$ is positive, moments cannot generally be computed so explicitly.  Taking advantage of the fact that the distributions, when all parameters except $\gvec \alpha$ are fixed, belong to the exponential family of distributions, it is easier to compute log-ratio covariances $\Cov [\log (X_i/X_j), \log(X_k/X_l)]$.  Moreover, the importance of log-ratios has been emphasized in the context of compositional data analysis \cite{aitchison}.

As a mathematical aside, the normalization constants of the distributions involve simplex integrals that relate to generalized hypergeometric functions, in particular the Carlson $R$-function, which is a rewriting of the Lauricella $F_D$-hypergeometric function.  In certain cases these provide generalizations of the complete homogeneous symmetric polynomials $h_n (X_1 , \ldots , X_K)$.  Taking a converse viewpoint, we find a natural way of generalizing such polynomials from non-negative integer $n$ to fractional $n$ through an integral representation on the simplex $S_K$.

In Section 2, we review some known probability distributions on the simplex, unifying them within a common framework.  Section 3 contains results of a more mathematical nature, for later use or to highlight connections to mathematical literature, in particular to symmetric polynomials.  In Section 4, we examine properties of the distributions, deriving moments and log-ratio covariances.  We conclude in Section 5.


\section{Probability distributions}


We first review some known probability distributions on the simplex and define our notation.  The probability simplex is denoted $S_K$, where the subscript $K$ denotes the dimension of the embedding space $\bb R^K$; the simplex itself is $(K - 1)$-dimensional.  This can be a more natural notation for applications, so the sum constraint $\sum_{i = 1}^K x_i = 1$ corresponds to a unit interval being divided into $K$ parts.

For $K$-dimensional vectors, we use the notation $\gvec \alpha = (\alpha_1 , \ldots , \alpha_K)$ and denote the sum of components $\alpha_+ = \sum_{i = 1}^K \alpha_i$.  The vectors that we consider have non-negative components, $\alpha_i \geq 0$, so $\alpha_+$ is also the $L^1$-norm of $\gvec \alpha$, i.e.\ $\alpha_+ = || \gvec \alpha ||_1$.  $\vec e_i$ is a unit-vector in the $i$-direction, i.e.\ with $j$-component $(\vec e_i)_j = \delta_{i j}$.  The $K$-dimensional vector whose entries are all 1 is $\vec 1 = (1 , \ldots, 1)$.  Given a vector $\gvec \beta$, whose components are $\beta_i$, the vector $\gvec \beta^\sigma$ is defined to have components $\beta_i^\sigma$, i.e.\ $\gvec \beta^\sigma = (\beta_1^\sigma , \ldots , \beta_K^\sigma)$.  We use bars to denote the operation of normalizing a vector to give a unit-vector, e.g.\ $\overline{\vec 1} = \vec 1 /K$ and $\overline{\gvec \beta ^{-\sigma}} = \gvec \beta^{-\sigma}/(\beta^{-\sigma})_+$.


\subsection{Dirichlet, Schl\"{o}milch and generalized beta distributions}



\subsubsection{Dirichlet distribution}


The $(K - 1)$-variate Dirichlet distribution is parameterized by a $K$-dimensional vector $\gvec \alpha = (\alpha_1 , \ldots , \alpha_K)$, with $\alpha_ i > 0$.  The probability density function, defined on the probability simplex $S_K$, is
\be
f(\vec x) = \fr{\Gamma (\alpha_+)}{\prod_{i = 1}^K \Gamma(\alpha_i)} \prod_{i = 1}^K x_i^{\alpha_i - 1} .
\label{dirichletdist}
\ee
The Dirichlet integral \eq{dirichletint} ensures the normalization $\int_{S_K} \! \df^{K - 1} x \, f(\vec x) = 1$, which is why Wilks gave the Dirichlet distribution its name \cite{wilks}.  For a random variable $\vec Y = (Y_1 , \ldots , Y_K)$ with this Dirichlet distribution, we denote $\vec Y \sim D(\gvec \alpha)$.  The special case $D(\vec 1)$ is simply a uniform distribution on the simplex.  Samples from the Dirichlet distribution $\vec Y$ can be generated from samples of standard gamma distributions $Z_i$ through
\be
Y_i = \fr{Z_i}{\sum_{j = 1}^K Z_j} ,
\label{dirichletgamma}
\ee
where $Z_i \sim \textrm{Gamma}(\alpha_i, 1)$ and $\alpha_i$ are shape parameters.

The mean, variance, and covariance for $i \neq j$ are
\begin{align}
\bb E (Y_i) & = \fr{\alpha_i}{\alpha_+} , & \Var(Y_i) & = \fr{\alpha_i (\alpha_+ - \alpha_i)}{\alpha_+^2 (\alpha_+ + 1)}, & \Cov (Y_i , Y_j) = - \fr{\alpha_i \alpha_j}{\alpha_+^2 (\alpha_+ + 1)} .
\end{align}
The $K$ independent parameters of $\gvec \alpha$ can be considered as corresponding to $K - 1$ independent means and an overall scaling of the covariance.  Note that $\Cov (Y_i , Y_j)$ is always negative, so the distribution is not appropriate for modelling data that has negative covariances.

For background on the Dirichlet distribution, including basic properties, applications and generalizations, see e.g.\ \cite{wilks, kobajo, ngtita}; a review of its history is \cite{gupric}.


\subsubsection{Schl\"{o}milch distribution}
\label{schlodistsection}


A more general family of distributions defined on $S_K$ is given by the probability density function
\be
f(\vec x) = \fr{\Gamma (\alpha_+) \prod_{j = 1}^K \beta_j^{\alpha_j}}{\prod_{i = 1}^K \Gamma(\alpha_i)} \fr{\prod_{i = 1}^K x_i^{\alpha_i - 1}}{(\sum_{j = 1}^K \beta_j x_j)^{\alpha_+}} .
\label{schlomilchdist}
\ee
Compared to the Dirichlet distribution, there are double the number of parameters, with an additional $K$-dimensional parameter vector $\gvec \beta = (\beta_1 , \ldots , \beta_K)$, with $\beta_i > 0$.  The distribution is invariant under rescaling $\gvec \beta \rightarrow \lambda \gvec \beta$, for $\lambda > 0$, so there are $2 K - 1$ independent parameters.  For identifiability when performing statistical inference, we may assume that $\sum_{i = 1}^K \beta_i = 1$, which implies that $\beta_i < 1$.  If all $\beta_i$ are equal, then the Dirichlet distribution is recovered.

In statistical literature, the distribution \eq{schlomilchdist} seems to have first appeared in work of Dickey \cite{dickey}, obtained from the Dirichlet distribution through the transformation \eq{transform} and attributed to a Savage communication.  Some other early appearances include the generalized Dirichlet distribution (G3D) of Chen and Novick \cite{chenov} and the scaled Dirichlet distribution of Aitchison \cite{aitchison}.  A more recent appearance is the gamma normalized infinitely divisible (NID) distribution \cite{fahapr, manben}, although one should note that the generalized gamma distribution, from which it is derived, is infinitely divisible for only certain ranges of its parameters.

A further parameter $\tau > 0$ may be introduced, giving a distribution with probability density function
\be
f(\vec x) = \fr{\tau^{K - 1} \Gamma (\alpha_+)}{(\sum_{j = 1}^K \beta_j x_j^\tau)^{\alpha_+}} \prod_{i = 1}^K \fr{\beta_i^{\alpha_i} x_i^{\tau \alpha_i - 1}}{\Gamma(\alpha_i)} .
\label{Sdist}
\ee
This distribution was obtained by Craiu and Craiu \cite{craiu2} (see also Chapter 49.1 of \cite{kobajo}), by transformation of generalized gamma distributions.  Equivalently, the transformation can be expressed in terms of standard gamma distributions as
\be
X_i = \fr{(Z_i/\beta_i)^{1/\tau}}{\sum_{j = 1}^K (Z_j/\beta_j)^{1/\tau}} , 
\label{transformXZ}
\ee
where $Z_i \sim \textrm{Gamma}(\alpha_i, 1)$ and $\alpha_i$ are shape parameters.

The case $\tau = 1$ corresponds to the coordinate transformation \eq{transform} and the distribution \eq{schlomilchdist}.  The distribution \eq{Sdist} has also been called the shifted-scaled Dirichlet distribution \cite{momapaeg} and the simplicial generalized beta distribution \cite{graf}.  The special case with $\gvec \beta = \overline{\vec 1}$ has also been called the generalized Dirichlet distribution \cite{dnposc}.

For the probability density function \eq{Sdist}, we shall denote the corresponding distribution by $S(\gvec \alpha, \gvec \beta, \tau)$.  The $S$ stands for Schl\"{o}milch, Savage, scaled, shifted, or simplicial, according to taste.  I shall refer to the $S$-distribution as the Schl\"{o}milch distribution, in order to complement the distribution named after Dirichlet.

Recalling the expression \eq{dirichletgamma}, which relates a Dirichlet random variable to gamma random variables, we see that there is a transformation relating a Dirichlet random variable $\vec Y \sim D(\gvec \alpha)$ and a Schl\"{o}milch random variable $\vec X \sim S(\gvec \alpha, \gvec \beta, \tau)$,
\begin{align}
Y_i & = \fr{\beta_i X_i^\tau}{\sum_{j = 1}^K \beta_j X_k^\tau} , & X_i & = \fr{(Y_i/\beta_i)^{1/\tau}}{\sum_{j = 1}^K (Y_j/\beta_j)^{1/\tau}} .
\label{duality}
\end{align}
In general, this transformation is trivial if and only if both $\tau = 1$ and all $\beta_i$ are equal.


\subsubsection{Univariate distributions: generalized beta}


In the univariate case, $K = 2$, the Dirichlet distribution is simply the beta distribution.  The three-parameter generalized beta distribution (G3B) of Libby and Novick \cite{libnov} matches the univariate Schl\"{o}milch distribution with $\tau = 1$, i.e.\ \eq{schlomilchdist}; later rediscoveries include \cite{seshadri, boshga, alpawo}.

More generally, the four-parameter generalized beta distribution (G4B) of Chen and Novick \cite{chenov} has probability density function of the form
\be
f(x) = \fr{\Gamma(\alpha_1 + \alpha_2)}{\Gamma(\alpha_1) \Gamma(\alpha_2) {_2} F {_1} (\alpha_1, \kappa; \alpha_+; 1 - \lambda)} \fr{x^{\alpha_1 - 1} (1 - x)^{\alpha_2 - 1}}{(\lambda x + 1 - x)^\kappa} ,
\label{g4b}
\ee
where $\alpha_i > 0$, $\kappa > 0$ and $\lambda > 0$.  Choosing $\kappa = \alpha_1 + \alpha_2$ recovers the G3B distribution.  The hypergeometric function ${_2} F {_1} (\alpha_1, \kappa; \alpha_+; 1 - \lambda)$ can be defined by demanding the normalization condition $\int_0^1 \! \df x \, f(x) = 1$; this corresponds to the Euler integral representation
\be
{_2}F{_1}(a, b; c; x) = \fr{\Gamma(c)}{\Gamma(b) \Gamma(c - b)} \int_0^1 \! \df t \, t^{b - 1} (1 - t)^{c - b - 1} (1 - x t)^{-a} ,
\ee
defined for $\Real c > \Real b > 0$ and then analytically continued.  Other early appearances of the distribution include the hypergeometric family of Chamayou and Letac \cite{chalet}, and the Gauss hypergeometric distribution of Armero and Bayeri \cite{armbay}.

In the case that $\kappa = \alpha_+ + n$, for a non-negative integer $n$, the G4B distribution is a finite mixture of Schl\"{o}milch distributions.  To see this, multiply the probability density function \eq{g4b} by $1 = (x + 1 - x)^n$ and use a binomial expansion, giving
\be
f(x) = \sum_{k = 0}^n w_k f_k(x) ,
\ee
where the weights are
\be
w_k = \binom{n}{k} \fr{(\alpha_1)_k (\alpha_2)_{n - k}}{(\alpha_+)_n \lambda^{\alpha_1 + k} {_2} F {_1} (\alpha_1, \alpha_+ + n; \alpha_+; 1 - \lambda)} ,
\ee
the Pochhammer symbol denotes the rising factorial, i.e.\ $(x)_n = x (x + 1) \ldots (x + n - 1)$, and
\be
f_k(x) = \fr{\Gamma(\alpha_+ + n) \lambda^{\alpha_1 + k}}{\Gamma(\alpha_1 + k) \Gamma(\alpha_2+ n - k)} \fr{x^{\alpha_1 + k - 1} (1 - x)^{\alpha_2 + n - k - 1}}{(\lambda x + 1 - x)^{\alpha_+ + n}}
\ee
is the probability density function of an $S(\alpha_1 + k, \alpha_2 + n - k, \frac{\lambda}{\lambda + 1}, \fr{1}{\lambda + 1}, 1)$ distribution.

Note that, because $\sum_{k = 0}^n w_k = 1$, the hypergeometric function can be expressed as a finite sum as
\be
{_2} F {_1} (\alpha_1, \alpha_+ + n; \alpha_+; 1 - \lambda) = \sum_{k = 0}^n \binom{n}{k} \fr{(\alpha_1)_k (\alpha_2)_{n - k}}{(\alpha_+)_n \lambda^{\alpha_1 + k}} .
\ee
This is equivalent to a known relation for the hypergeometric function,
\be
{_2} F {_1} (\alpha_1, \alpha_+ + n; \alpha_+; 1 - \lambda) = \fr{(\alpha_2)_n \, {_2}F{_1}(\alpha_1, -n; 1 - n - \alpha_2; \lambda^{-1})}{(\alpha_+)_n \lambda^{\alpha_1} } ,
\ee
by using the definition of the hypergeometric function on the right as a power series that converges for $|\lambda^{-1}| < 1$, analytically continued.  Although a known result, I have been unable to find a direct derivation in easily accessible literature.  It may be derived from known properties of the hypergeometric function, using the Euler transformation
\be
{_2}F{_1}(a, b; c; x) = (1 - x)^{c - a - b} {_2}F{_1}(c - a, c - b; c; x) ,
\ee
and taking the $a \rightarrow -n$ limit, with $b = \alpha_2$, $c = \alpha_1 + \alpha_2$ and $x = 1 - \lambda$, in the relation (see e.g.\ Chapter 1.8 of \cite{slater})
\begin{align}
{_2}F{_1}(a, b; c; x) & = \fr{\Gamma(b - a) \Gamma(c)}{\Gamma(b) \Gamma(c - a)} (1 - x)^{-a} {_2}F{_1}(a, c - b; a - b + 1; 1/(1 - x)) \nnr
& \qquad + \fr{\Gamma(a - b) \Gamma(c)}{\Gamma(a) \Gamma(c - b)} (1 - x)^{-b} {_2}F{_1}(b, c - a; b - a + 1; 1/(1 - x)) ,
\end{align}
which is a linear relation between three of Kummer's 24 solutions of the hypergeometric equation, a linear, second-order differential equation.

A different generalization of the G3B distribution is given by the $GB(a, b, c, \alpha_1, \alpha_2)$ distribution of McDonald and Xu \cite{mcdxu}, which has support on $(0, b/(1 - c)^{1/a})$, where $b > 0$ and $0 \leq c \leq 1$.  Choosing $a > 0$ and $b = (1 - c)^{1/a}$, with $c \neq 1$, we have a 4-parameter family of distributions with support on $(0, 1)$, with probability density function
\be
f(x) = \fr{\Gamma(\alpha_1 + \alpha_2) a x^{a \alpha_1 - 1} (1 - x^a)^{\alpha_2 - 1}}{\Gamma(\alpha_1) \Gamma(\alpha_2) (1 - c)^{\alpha_1} [1 + c x^a/(1 - c)]^{\alpha_1 + \alpha_2}} .
\ee
This is the univariate case, $K = 2$, of the Schl\"{o}milch distribution \eq{Sdist}.  The case $a = 1$ reduces to the G3B distribution.


\subsection{Mixture distributions}



\subsubsection{Normalization constant for Dirichlet mixture}


Consider a special Dirichlet mixture distribution, which we denote $DM(\gvec \alpha, \gvec \gamma, n, \sigma)$, with probability density function of the form
\be
f(\vec x) = \fr{1}{I_n^\sigma (\gvec \alpha, \gvec \gamma)} \bigg( \sum_{k = 1}^K \gamma_k x_k^\sigma \bigg) ^n \prod_{i = 1}^K x_i^{\alpha_i - 1} ,
\label{DM}
\ee
where $n$ is a positive integer, $\sigma > 0$ and, for identifiability, we restrict the parameter vector $\gvec \gamma$ to have unit-norm, i.e.\ $\gamma_+ = 1$.  We define the function
\be
I_n^\sigma (\gvec \alpha, \gvec \gamma) = \int_{S_K} \! \df^{K - 1} y \, \bigg( \sum_{k = 1}^K \gamma_k y_k^\sigma \bigg) ^n \prod_{i = 1}^K y_i^{\alpha_i - 1} ,
\label{Ins}
\ee
where the arguments satisfy $\alpha_i > 0$, $\gamma_i > 0$, $\sigma > 0$, which plays the role of a normalization constant for the Dirichlet mixture distribution.  For most of our study, it suffices to take $n$ to be a non-negative integer, although the function can be extended to real $n$ and, with a choice of branch cut, to complex $n$.  The function $I_n^\sigma(\gvec \alpha, \gvec \gamma)$ satisfies the scaling $I_n^\sigma(\gvec \alpha, \lambda \gvec \gamma) = \lambda^n I_n^\sigma(\gvec \alpha, \gvec \gamma)$, for $\lambda > 0$, implying that $I_n^\sigma(\gvec \alpha, \gvec \gamma) = \gamma_+^n I_n^\sigma(\gvec \alpha, \overline{\gvec \gamma})$.

The Dirichlet mixture distribution $DM(\gvec \alpha, \gvec \gamma, n, \sigma)$ can be explicitly expressed as a mixture of Dirichlet distributions by decomposing its probability density function \eq{DM} as
\be
f(\vec x) = \sum_{n_+ = n} w_{\vec n} f_{\vec n} (\vec x) , 
\ee
where $f_{\vec n} (\vec x)$ is the probability density function of a Dirichlet distribution $D(\gvec \alpha + \sigma \vec n)$ obtained from \eq{dirichletdist}.  The summation is over non-negative integers $n_1 , \ldots , n_K$ such that $\sum_{i = 1}^K n_i = n$.  The weights are
\be
w_{\vec n} = \fr{1}{I_n^\sigma (\gvec \alpha, \gvec \gamma) \Gamma(\alpha_+ + \sigma n)} \binom{n}{\vec n} \prod_{i = 1}^K \Gamma(\alpha_i + \sigma n_i) \gamma_i^{n_i} ,
\label{weights}
\ee
using the notation for a multinomial coefficient
\be
\binom{n}{\vec m} = \fr{n!}{\prod_{i = 1}^K m_i!} ,
\ee
where $\sum_{i = 1}^K m_i = n$.  The overall normalization of the Dirichlet mixture distribution implies that $\sum_{\vec n} w_{\vec n} = 1$, and so we deduce an expression for $I_n^\sigma (\gvec \alpha, \gvec \gamma)$ as a finite sum,
\be
I_n^\sigma(\gvec \alpha, \gvec \gamma) = \fr{1}{\Gamma(\alpha_+ + \sigma n)} \sum_{n_+ = n} \binom{n}{\vec n} \prod_{i = 1}^K \Gamma(\alpha_i + \sigma n_i) \gamma_i^{n_i} ,
\label{Inmultinomial}
\ee
which relies on $n$ being a non-negative integer.  In some cases, which we consider later, this sum may be further simplified.


\subsubsection{Jacobian computation}


We perform a coordinate change from $\vec y$ to $\vec x$ coordinates by defining
\be
x_i = \fr{(y_i/\beta_i)^{1/\tau}}{\sum_{j = 1}^K (y_j/\beta_j)^{1/\tau}} ,
\label{ytox}
\ee
for $\tau > 0$, corresponding to the transformation of random variables \eq{duality}.  It is clear that $\sum_{i = 1}^K x_i = 1$, so $\vec x$ are coordinates on the simplex $S_K$.  The inverse coordinate transformation is
\be
y_i = \fr{\beta_i x_i^\tau}{\sum_{j = 1}^K \beta_j x_j^\tau} .
\label{xtoy}
\ee
By computing the Jacobian, we obtain the transformation of the measure
\be
\fr{\df^{K - 1} y}{\prod_{i = 1}^K y_i} = \tau^{K - 1} \fr{\df^{K - 1} x}{\prod_{i = 1}^K x_i} .
\label{measure}
\ee

It is a routine calculation to use the explicit coordinate transformations of \eq{ytox} or \eq{xtoy} and compute partial derivatives $\pd y_i/\pd x_j$ or $\pd x_i/\pd y_j$ to find the Jacobian. However, a more implicit expression of the coordinate change gives a clearer understanding of the origin of the transformation of the measure \eq{measure}.  From \eq{ytox} or \eq{xtoy}, we have the differential relation 
\be
\df \log \bigg( \fr{y_i}{y_{K - 1}} \bigg)  = \tau \, \df \log \bigg( \fr{x_i}{x_{K - 1}} \bigg) .
\label{diffrel}
\ee
Remembering that $y_K$ is expressed in terms of the independent coordinates as $y_K = 1 - \sum_{i = 1}^{K - 1} y_i$, we take exterior products to compute the volume form (Chapter 2 of \cite{muirhead} provides an introduction to differential forms in a statistical context)
\begin{align}
& \df \log \bigg( \fr{y_1}{y_{K - 1}} \bigg) \wedge \df \log \bigg( \fr{y_2}{y_{K - 1}} \bigg) \wedge \ldots \wedge \df \log \bigg( \fr{y_{K - 2}}{y_{K - 1}} \bigg) \wedge \df \log \bigg( \fr{y_K}{y_{K - 1}} \bigg) \nnr
& = \bigg( \fr{\df y_1}{y_1} - \fr{\df y_{K - 1}}{y_{K - 1}} \bigg) \wedge \bigg( \fr{\df y_2}{y_2} - \fr{\df y_{K - 1}}{y_{K - 1}} \bigg) \wedge \ldots \wedge \bigg( \fr{\df y_{K - 2}}{y_{K - 2}} - \fr{\df y_{K - 1}}{y_{K - 1}} \bigg) \wedge \bigg( - \sum_{i = 1}^{K - 1} \fr{\df y_i}{y_K} - \fr{\df y_{K - 1}}{y_{K - 1}} \bigg) \nnr
& = \fr{\df y_1 \wedge \df y_2 \wedge \ldots \wedge \df y_{K - 1}}{\prod_{i = 1}^{K - 1} y_i} \bigg( - 1 - \fr{\sum_{i = 1}^{K - 1} y_i}{y_K} \bigg) \nnr
& = - \fr{\df y_1 \wedge \df y_2 \wedge \ldots \wedge \df y_{K - 1}}{\prod_{i = 1}^K y_i} .
\end{align}
From the same calculation with $\vec x$ instead of $\vec y$, and the differential relation \eq{diffrel}, we obtain
\be
\fr{\df y_1 \wedge \df y_2 \wedge \ldots \wedge \df y_{K - 1}}{\prod_{i = 1}^K y_i} = \tau^{K - 1} \fr{\df x_1 \wedge \df x_2 \wedge \ldots \wedge \df x_{K - 1}}{\prod_{i = 1}^K x_i} ,
\ee
which is the transformation of the measure \eq{measure} equivalently expressed in terms of volume forms.

Note that the computation works in the same way for the opposite sign of $\tau$, by replacing $\tau \rightarrow - \tau$, giving an extra factor of $(-1)^{K - 1}$.  When computing an integral over the simplex, this factor is cancelled out by consistently ordering the limits of the integral.


\subsubsection{Schl\"{o}milch mixture distribution}


The coordinate transformation leads to the special Schl\"{o}milch mixture distribution, which we denote by 
\be
\vec X \sim SM(\gvec \alpha, \gvec \beta, \gvec \gamma, n, \sigma, \tau) .
\ee
Using the transformation of the measure given by \eq{measure}, we obtain the probability density function
\be
f(\vec x) = \fr{\tau^{K - 1} ( \prod_{l = 1}^K \beta_l^{\alpha_l} ) ( \sum_{k = 1}^K \beta_k^\sigma \gamma_k x_k^{\sigma \tau} ) ^n \prod_{i = 1}^K x_i^{\tau \alpha_i - 1}}{I_n^\sigma (\gvec \alpha, \gvec \gamma) (\sum_{j = 1}^K \beta_j x_j^\tau)^{\alpha_+ + \sigma n}} ,
\label{SM}
\ee
where $n$ is a positive integer, $\sigma > 0$, $\tau > 0$.  For identifiability, we restrict the parameter vectors $\gvec \beta$ and $\gvec \gamma$ to have unit-norm, i.e.\ $\beta_+ = \gamma_+ = 1$, so there are $3K$ independent continuous parameters, in addition to the non-negative integer parameter $n$.

The Schl\"{o}milch mixture distribution can be explicitly expressed as a mixture of Schl\"{o}milch distributions by decomposing its probability density function \eq{SM} as
\be
f(\vec x) = \sum_{n_+ = n} w_{\vec n} f_{\vec n} (\vec x) , 
\ee
where $f_{\vec n} (\vec x)$ is the probability density function of a Schl\"{o}milch distribution $S(\gvec \alpha + \sigma \vec n, \gvec \beta, \tau)$ obtained from \eq{schlomilchdist}.  The weights $w_{\vec n}$ are \eq{weights}, the same as for the original Dirichlet mixture distribution.


\subsubsection{Relations between distributions}

\begin{figure}
\centering
\includegraphics{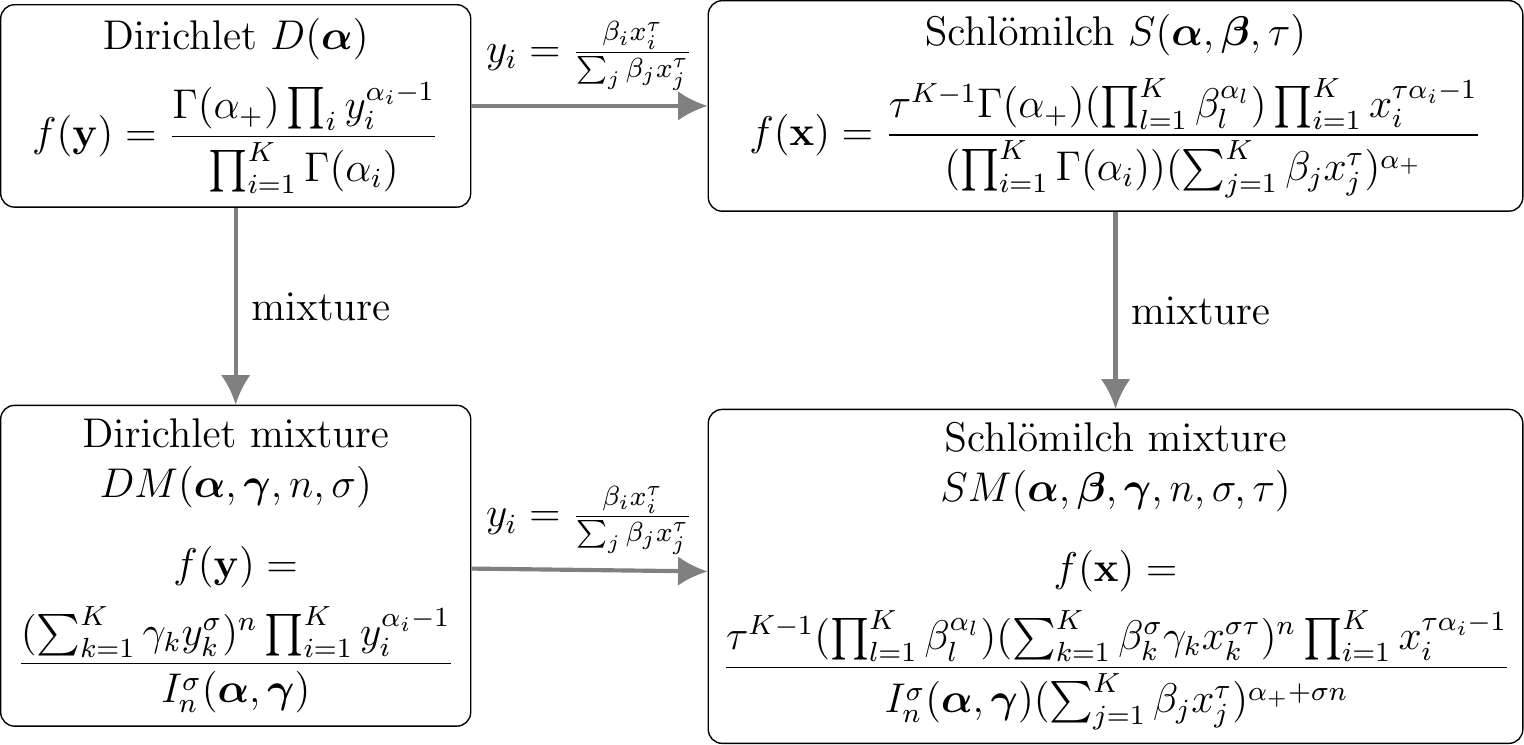}
\caption{Two routes from the Dirichlet distribution to the Schl\"{o}milch mixture distribution}
\label{4d}
\end{figure}

The connections between the four key families of distributions are illustrated in Figure \ref{4d}.  To generate samples from the Schl\"{o}milch mixture distribution, one can generate samples from the Dirichlet distribution, followed by taking mixtures and transforming.

Conversely, one may regard the Schl\"{o}milch mixture distribution as providing a unified framework for several distributions that correspond to special choices of parameters, as seen from the probability density functions in Figure \ref{4d}, remembering the constraint $\sum_{i = 1}^K x_i = 1$.  For each positive integer $n$, the Schl\"{o}milch mixture distribution $SM(\gvec \alpha, \gvec \beta, \gvec \gamma, n, \sigma, \tau)$ contains, as special cases, the Dirichlet mixture distribution $DM(\gvec \alpha, \gvec \gamma, n, \sigma)$, the Schl\"{o}milch distribution $S(\gvec \alpha, \gvec \beta, \tau)$ and the Dirichlet distribution $D(\gvec \alpha)$.

More generally, there are three special cases in which the probability density function of the Schl\"{o}milch mixture distribution \eq{SM} simplifies.  The first case is when the denominator factor $\sum_{j = 1}^K \beta_j x_j^\tau$ is constant, corresponding to the Dirichlet mixture distribution.  The second case is when the numerator factor $\sum_{k = 1}^K \beta_k^\sigma \gamma_k x_k^{\sigma \tau}$ and the denominator factor $\sum_{j = 1}^K \beta_j x_j^\tau$ are proportional, leaving only the denominator factor with a reduced exponent, which is the Schl\"{o}milch distribution.  The third case is when the numerator factor $\sum_{j = 1}^K \beta_j^\sigma \gamma_j x_j^{\sigma \tau}$ is constant, by taking $\gamma_i = \beta_i^{-\sigma}/(\beta^{-\sigma})_+$ and $\sigma \tau = 1$, giving the probability density function
\be
f(\vec x) = \fr{( \prod_{l = 1}^K \beta_l^{\alpha_l} )  \prod_{i = 1}^K x_i^{\alpha_i/\sigma - 1}}{\sigma^{K - 1}  I_n^{\sigma} (\gvec \alpha, \gvec \beta^{-\sigma}) (\sum_{j = 1}^K \beta_j x_j^{1/\sigma})^{\alpha_+ + \sigma n}} .
\label{case3}
\ee
In the univariate $K = 2$ case with $\sigma = 1$ and $n$ generalized to take continuous values, this corresponds to the G4B distribution of \eq{g4b}.

\begin{figure}
\centering
\includegraphics{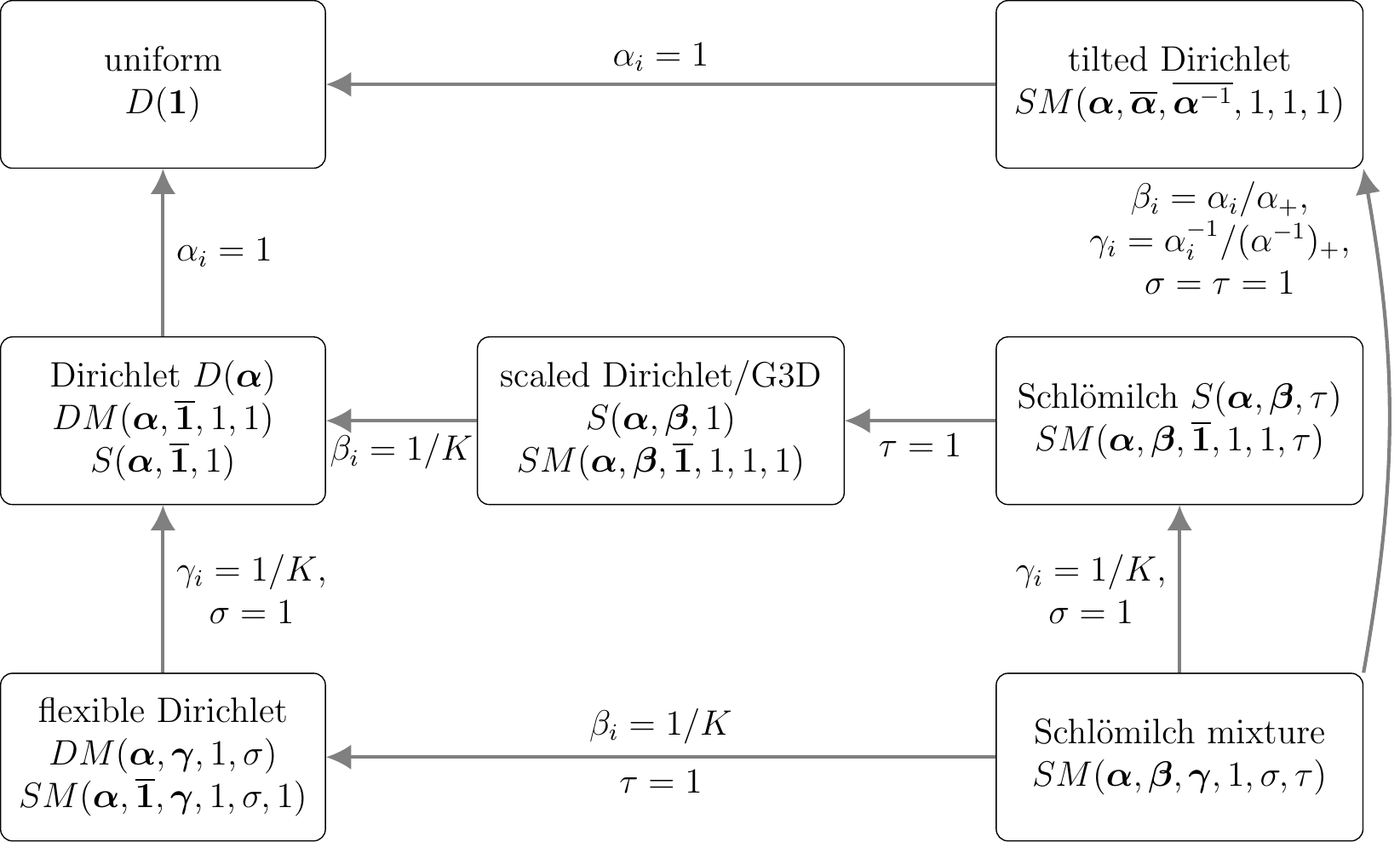}
\caption{Relations between $n = 1$ distributions}
\label{4dn1}
\end{figure}

The most basic $n = 1$ case illustrates how the Dirichlet, Dirichlet mixture and Schl\"{o}milch distributions arise as special cases of the Schl\"{o}milch mixture distribution, as shown in Figure \ref{4dn1}.  In this $n = 1$ case, the Dirichlet mixture distribution, which is a general mixture of $D(\gvec \alpha + \sigma \vec e_i)$ distributions, has been called the flexible Dirichlet distribution by Ongaro, Migliorati and Monti \cite{onmimo} (\cite{ngtita} prefers to call it the mixed Dirichlet distribution).  The flexible Dirichlet distribution has been further studied in \cite{ongmig, mionmo} and, unlike the Dirichlet distribution, allows for multi-modality.  Figure \ref{4dn1} explicitly includes the uniform distribution $D(\vec 1)$ and the scaled Dirichlet or G3D distribution, as discussed in Section \ref{schlodistsection}, which arise as special cases of the Schl\"{o}milch mixture distribution for any choice of $n$.

In addition, Figure \ref{4dn1} includes $SM(\gvec \alpha, \overline{\gvec \alpha}, \overline{\gvec \alpha^{-1}}, 1, 1, 1)$, i.e.\ \eq{case3} with $\gvec \beta = \overline{\gvec \alpha}$, $n = 1$ and $\tau = 1$, for which $\bb E X_i = 1/K$ (so $\bb E \sum_{i = 1}^K X_i = 1$).  This distribution has been used in the context of multivariate extreme value theory by Coles and Tawn \cite{coltaw}, and has subsequently been referred to as the tilted Dirichlet distribution.  The normalization involves $I_1^1 (\gvec \alpha, \gvec \alpha^{-1}) = K \prod_{i = 1}^K \Gamma(\alpha_i)/\Gamma (1 + \alpha_+)$, as known to Schl\"{o}milch in \eq{schloint1}, and so the probability density function is
\be
f(\vec x) = \fr{\Gamma(1 + \alpha_+)}{K} \bigg( \prod_{l = 1}^K \fr{\alpha_l^{\alpha_l}}{\Gamma (\alpha_l)} \bigg) \fr{\prod_{i = 1}^K x_i^{\alpha_i - 1}}{(\sum_{j = 1}^K \alpha_j x_j)^{\alpha_+ + 1}} .
\ee

The $n = 2$ case has received less attention in the literature.  However, the recently considered double flexible Dirichlet distribution \cite{asmion} is a general mixture of $D(\gvec \alpha + \sigma(\vec e_i + \vec e_j))$ distributions, involving $\binom{K}{2} - 1$ independent mixture weights.  This includes the more restricted $DM(\gvec \alpha, \gvec \gamma, 2, \sigma)$ distribution, which involves only $K - 1$ independent mixture weights.


\subsection{Inverse Schl\"{o}milch distribution}


Consider the distribution constructed by taking
\be
X_i = \fr{(\beta_i/Z_i)^{1/\tau}}{\sum_{j = 1}^K (\beta_j/Z_j)^{1/\tau}} , 
\label{transformXZ2}
\ee
where $Z_i \sim \textrm{Gamma}(\alpha_i, 1)$ and $\alpha_i$ are shape parameters.  Again, the components of $\gvec \beta$ are positive, $\beta_i > 0$, and we shall assume that $\sum_{i = 1}^K \beta_i = 1$ for identifiability.  This is simply the transformation used previously in \eq{transformXZ}, with the replacement $\tau \rightarrow - \tau$; however, it is convenient to state this transformation separately, continuing to take $\tau > 0$, for comparison with literature.  Using the transformation of the measure given by \eq{measure}, but with $\tau \rightarrow - \tau$, we obtain the probability density function
\be
f(\vec x) = \fr{\tau^{K - 1} \Gamma (\alpha_+)}{(\sum_{j = 1}^K \beta_j/x_j^\tau)^{\alpha_+}} \prod_{i = 1}^K \fr{\beta_i^{\alpha_i}}{\Gamma(\alpha_i) x_i^{\tau \alpha_i + 1}} .
\ee
We shall call this an inverse Schl\"{o}milch distribution $IS(\gvec \alpha, \gvec \beta, \tau)$.

The special case $IS(\vec 1 , \gvec \beta, \tau)$ has applications in neural networks, where it has been called the Concrete distribution \cite{mamnte} (``Concrete'' is a portmanteau of ``continuous'' and ``discrete'') and the Gumbel-softmax distribution \cite{jagupo}.  If $\alpha_i = 1$, then $Z_i \sim \textrm{Exponential}(1)$ is a unit exponential distribution and $W_i = - \log Z_i \sim \textrm{Gumbel}(0, 1)$ is a standard Gumbel distribution.  Then the transformation \eq{transformXZ2} takes the form
\be
X_i = \textrm{softmax} [(\vec W + \log \gvec \beta)/\tau]_i , 
\ee
where $\log \gvec \beta = (\log \beta_1 , \ldots , \log \beta_K)$ and the softmax function, given by
\be
\textrm{softmax}(\vec x)_i = \fr{\exp (x_i)}{\sum_{j = 1}^K \exp (x_j)} ,
\ee
is a mapping from $\bb R^K$ to $S_K \subset \bb R^K$.  In this context, the parameter $\tau$ is called a temperature, and the $\tau \rightarrow 0$ limit tends to a categorical distribution.  The vector $\gvec \beta$ is related to the corresponding probabilities of a categorical distribution; the constraint $\sum_{i = 1}^K \beta_i = 1$ ensures that they are normalized as probabilities.  The more general inverse Schl\"{o}milch distribution $IS(\gvec \alpha, \gvec \beta, \tau)$ may therefore be considered as a Dirichlet--Concrete or Dirichet--Gumbel-softmax distribution.

The special case $IS(\gvec \alpha, \gvec \beta, 1)$ has been considered by del Castillo \cite{delcas}, through the same construction but expressed in terms of chi-squared distributions instead of gamma distributions.  Within the del Castillo distribution, the case $IS(\tf{1}{2} \vec 1, \gvec \beta, 1)$ is an example from the ratio-stable distributions considered by Carlton \cite{carltonphd, carlton}, parameterized by $\gvec \beta$ and one extra parameter, which are related to a Dirichlet process.

In the univariate $K = 2$ case, the inverse Schl\"{o}milch distribution $IS(\gvec \alpha, \gvec \beta, \tau)$ has probability density function
\begin{align}
f(x) & = \fr{\tau \Gamma(\alpha_1 + \alpha_2) \beta_1^{\alpha_1} \beta_2^{\alpha_2}}{\Gamma(\alpha_1) \Gamma(\alpha_2) x^{\tau \alpha_1 + 1} (1 - x)^{\tau \alpha_2 + 1} [\beta_1/x^\tau + \beta_2/(1 - x)^\tau]^{\alpha_+}} \nnr
& = \fr{\tau \Gamma(\alpha_1 + \alpha_2) \beta_1^{\alpha_1} \beta_2^{\alpha_2} x^{\tau \alpha_2 - 1} (1 - x)^{\tau \alpha_1 - 1}}{\Gamma(\alpha_1) \Gamma(\alpha_2) [\beta_2 x^\tau + \beta_1 (1 - x)^\tau]^{\alpha_+}} .
\end{align}
This is the same form as a univariate Schl\"{o}milch distribution, since $IS(\alpha_1, \alpha_2, \beta_1, \beta_2, \tau) = S(\alpha_2, \alpha_1, \beta_2, \beta_1, \tau)$.

There does not appear to be a construction of special mixture distributions for the inverse Schl\"{o}milch distribution that is analogous to the construction for the Schl\"{o}milch distribution.


\section{Normalization constant and symmetric polynomials}


In this section, we collect some useful results concerning the function $I_n^\sigma (\gvec \alpha, \gvec \gamma)$ that appears in the normalization of the distributions.  We also study in detail the special case of $\sigma = 1$, for which the function $I_n^1 (\gvec \alpha, \gvec \gamma)$ is, up to a $\gvec \gamma$-independent scaling, the Carlson $R$-function, and can be expressed in terms of deformations of standard symmetric polynomials.


\subsection{Basic relations}



\subsubsection{Dual integral representations}


From the normalization of the Schl\"{o}milch mixture distribution \eq{SM}, we see that the function $I_n^\sigma (\gvec \alpha, \gvec \gamma)$ has the equivalent integral representations
\be
I_n^\sigma (\gvec \alpha, \gvec \gamma) = \tau^{K - 1} ( {\textstyle \prod_{l = 1}^K \beta_l^{\alpha_l} } ) \int_{S_K} \! \df^{K - 1} x \, \fr{ ( \sum_{k = 1}^K \beta_k^\sigma \gamma_k x_k^{\sigma \tau} ) ^n \prod_{i = 1}^K x_i^{\tau \alpha_i - 1}}{(\sum_{j = 1}^K \beta_j x_j^\tau)^{\alpha_+ + \sigma n}} .
\ee
The redundancy of $\gvec \beta$ and $\tau$ here gives the freedom to make different choices of $\gvec \beta$ and $\tau$, which may provide more convenient integral representations of $I_n^\sigma (\gvec \alpha, \gvec \gamma)$ for different purposes.  The original integral representation \eq{Ins} is recovered by setting all $\beta_i$ equal and $\tau = 1$.  Another simple integral representation comes by setting $\beta_i = 1/\gamma_i^{1/\sigma}$ and $\tau = 1/\sigma$, giving
\be
I_n^\sigma (\gvec \alpha, \gvec \gamma) = \fr{1}{\sigma^{K - 1} \prod_{l = 1}^K \gamma_l^{\alpha_l/\sigma}} \int_{S_K} \! \df ^{K - 1} x \, \fr{\prod_{i = 1}^K x_i^{\alpha_i/\sigma - 1}}{(\sum_{j = 1}^K x_j^{1/\sigma}/\gamma_j^{1/\sigma})^{\alpha_+ + \sigma n}} .
\label{schloint2}
\ee
By comparison with the original integral representation for $I_n^\sigma (\gvec \alpha, \gvec \gamma)$ of \eq{Ins}, we obtain a duality relation
\be
I_n^\sigma (\gvec \alpha, \gvec \gamma) = \fr{I_{-(\alpha_+ + \sigma n)}^{1/\sigma} (\gvec \alpha/\sigma, \gvec \gamma^{-1/\sigma})}{\sigma^{K - 1} \prod_{i = 1}^K \gamma_i^{\alpha_i/\sigma}} .
\label{Idual}
\ee
The Pfaff transformation for the ${_2}F{_1}$ hypergeometric function,
\be
{_2}F{_1} (a, b; c; x) = (1 - x)^{-a} {_2}F{_1}(a, c - b; c; x/(1 - x)) ,
\ee
is a special case of the duality \eq{Idual} for $K = 2$ and $\sigma = 1$, which can be seen by using the Euler integral representation of the ${_2}F{_1}$ hypergeometric function.


\subsubsection{Recurrence relations}


By differentiating the integral representation \eq{schloint2}, we find that
\be
\bigg( - \fr{\pd}{\pd \gamma_k^{-1/\sigma}} \bigg) ^m \Big( {\textstyle (\prod_{l = 1}^K \gamma_l^{\alpha_l/\sigma})} I_n^\sigma (\gvec \alpha, \gvec \gamma) \Big) = \fr{(\alpha_+ + \sigma n)_m}{\sigma^{K - 1}} \int_{S_K} \! \df ^{K - 1} x \, \fr{x_k^{m/\sigma} \prod_{i = 1}^K x_i^{\alpha_i/\sigma - 1}}{(\sum_{j = 1}^K x_j^{1/\sigma}/\gamma_j^{1/\sigma})^{\alpha_+ + \sigma n + m}} .
\ee
In the case that $\sigma = m$ is a positive integer, this gives the differential recurrence relation
\be
{\textstyle (\prod_{l = 1}^K \gamma_l^{\alpha_l/m})} I_{n + 1}^m (\gvec \alpha, \gvec \gamma) = \fr{m^m}{(\alpha_+ + n m)_m} \sum_{k = 1}^K \bigg( \gamma_k^{1 + 1/m} \fr{\pd}{\pd \gamma_k} \bigg) ^m \Big( {\textstyle (\prod_{l = 1}^K \gamma_l^{\alpha_l/m})} I_n^m (\gvec \alpha, \gvec \gamma) \Big) ,
\label{diffrecurm}
\ee
a technique originally suggested by Schl\"{o}milch.

For general $\sigma$, the original integral representation \eq{Ins} immediately gives the algebraic recurrence relation
\be
I_{n + 1}^\sigma (\gvec \alpha, \gvec \gamma) = \sum_{i = 1}^K \gamma_i I_n^\sigma (\gvec \alpha + \sigma \vec e_i, \gvec \gamma) .
\label{algrecur}
\ee


\subsection{Low $n$ cases}
\label{lown}


For later use, it is convenient to record here some expressions involving $I_n^\sigma (\gvec \alpha, \gvec \gamma)$, as given by the finite sum \eq{Inmultinomial}, for low values of $n$.


\subsubsection{$n = 0$}


The Dirichlet integral \eq{dirichletint} gives
\be
I_0^\sigma (\gvec \alpha, \gvec \gamma) = \fr{\prod_{i = 1}^K \Gamma(\alpha_i)}{\Gamma (\alpha_+)} ,
\label{I0}
\ee
which is independent of $\sigma$ and $\gvec \gamma$.  Taking derivatives, we obtain
\be
\fr{\pd \log I_0^\sigma (\gvec \alpha, \gvec \gamma)}{\pd \alpha_i} = - \psi (\alpha_+) + \psi(\alpha_i)
\ee
and
\be
\fr{\pd^2 \log I_0^\sigma (\gvec \alpha, \gvec \gamma)}{\pd \alpha_i \, \pd \alpha_j} = - \psi' (\alpha_+) + \psi'(\alpha_i) \delta_{i j} ,
\label{d20}
\ee
where $\psi(x) = \Gamma'(x)/\Gamma(x) = (\log \Gamma (x))'$ is the digamma function, whose derivative $\psi'(x)$ is the trigamma function.


\subsubsection{$n = 1$}


The $n = 1$ case is relevant to several distributions studied in the literature, as indicated by Figure \ref{4dn1}.  In this case, we have
\be
I_1^\sigma (\gvec \alpha, \gvec \gamma) = \fr{\prod_{i = 1}^K \Gamma(\alpha_i)}{\Gamma(\alpha_+ + \sigma)} \sum_{k = 1}^K \gamma_k \fr{\Gamma(\alpha_k + \sigma)}{\Gamma(\alpha_k)} .
\ee
Taking derivatives, we obtain
\be
\fr{\pd \log I_1^\sigma(\gvec \alpha, \gvec \gamma)}{\pd \alpha_i}  = - \psi(\alpha_+ + \sigma) + \psi(\alpha_i) + p_i [\psi (\alpha_i + \sigma) - \psi(\alpha_i)] ,
\ee
where
\be
p_i = \fr{\gamma_i \Gamma(\alpha_i + \sigma)/\Gamma(\alpha_i)}{\sum_{k = 1}^K \gamma_k \Gamma(\alpha_k + \sigma)/\Gamma(\alpha_k)} 
\ee
have the interpretation of weights in the mixture distribution, as in \eq{weights}.  Further calculation gives
\be
\fr{\pd^2 \log I_1^\sigma(\gvec \alpha, \gvec \gamma)}{\pd \alpha_i \, \pd \alpha_j}  = - \psi'(\alpha_+ + \sigma) + \bigg( p_i \psi'(\alpha_i + \sigma) + (1 - p_i) \psi'(\alpha_i) + \fr{c_i^2}{p_i} \bigg) \delta_{i j} - c_i c_j ,
\label{d21}
\ee
where $c_i = p_i [\psi(\alpha_i + \sigma) - \psi(\alpha_i)]$.  

In the special case that additionally $\sigma = 1$, we have the further simplification
\be
I_1^1 (\gvec \alpha, \gvec \gamma) = \fr{\prod_{i = 1}^K \Gamma(\alpha_i)}{\Gamma (\alpha_+ + 1)} h_1 (\gvec \alpha, \gvec \gamma) ,
\label{I11}
\ee
where $h_1(\gvec \alpha, \gvec \gamma) = \sum_{i = 1}^K \alpha_i \gamma_i$.  The derivatives simplify as
\be
\fr{\pd \log I_1^1 (\gvec \alpha, \gvec \gamma)}{\pd \alpha_i} = - \psi (\alpha_+ + 1) + \psi(\alpha_i) + \fr{\gamma_i}{h_1 (\gvec \alpha, \gvec \gamma)}
\ee
and
\be
\fr{\pd^2 \log I_1^1 (\gvec \alpha, \gvec \gamma)}{\pd \alpha_i \, \pd \alpha_j} = - \psi' (\alpha_+ + 1) + \psi'(\alpha_i) \delta_{i j} - \fr{\gamma_i \gamma_j}{h_1 (\gvec \alpha, \gvec \gamma)^2} .
\label{d21b}
\ee
If we further set $\gvec \gamma = \vec 1$, then we recover the results for the Dirichlet integral, since $I_1^1(\gvec \alpha, \vec 1) = I_0^\sigma (\gvec \alpha, \gvec \gamma)$.


\subsubsection{$n = 2$}


For the $n = 2$ case, and further specializing to $\sigma = 1$, we have
\be
I_2^1 (\gvec \alpha, \gvec \gamma) = \fr{\prod_{i = 1}^K \Gamma(\alpha_i)}{\Gamma (\alpha_+ + 2)} [ p_1(\gvec \alpha, \gvec \gamma)^2 + p_2(\gvec \alpha, \gvec \gamma) ] = \fr{2 \prod_{i = 1}^K \Gamma(\alpha_i)}{\Gamma (\alpha_+ + 2)} h_2 (\gvec \alpha, \gvec \gamma) ,
\label{I21}
\ee
where $p_1(\gvec \alpha, \gvec \gamma) = \sum_{i = 1}^K \alpha_i \gamma_i$, $p_2(\gvec \alpha, \gvec \gamma) = \sum_{i = 1}^K \alpha_i \gamma_i^2$, and $h_2(\gvec \alpha, \gvec \gamma) = \tf{1}{2} [p_1(\gvec \alpha, \gvec \gamma)^2 + p_2 (\gvec \alpha, \gvec \gamma)]$. Taking derivatives, we obtain
\be
\fr{\pd \log I_2^1 (\gvec \alpha, \gvec \gamma)}{\pd \alpha_i} = - \psi (\alpha_+ + 2) + \psi(\alpha_i) + \fr{2 \gamma_i p_1 (\gvec \alpha, \gvec \gamma) + \gamma_i^2}{2 h_2 (\gvec \alpha, \gvec \gamma)} ,
\ee
and
\begin{align}
\fr{\pd^2 \log I_2^1 (\gvec \alpha, \gvec \gamma)}{\pd \alpha_i \, \pd \alpha_j} & = - \psi' (\alpha_+ + 2) + \psi'(\alpha_i) \delta_{i j} \nnr
& \qquad + \fr{\gamma_i \gamma_j [2 p_2(\gvec \alpha, \gvec \gamma) - 2 p_1(\gvec \alpha, \gvec \gamma)^2 - 2 (\gamma_i + \gamma_j) p_1(\gvec \alpha, \gvec \gamma) - \gamma_i \gamma_j]}{4 h_2(\gvec \alpha, \gvec \gamma)^2} .
\label{d22}
\end{align}
Setting $\gvec \gamma = \vec 1$ again recovers results for the Dirichlet integral, since $I_n^1 (\gvec \alpha, \vec 1) = I_0^\sigma(\gvec \alpha, \gvec \gamma)$ for any positive integer $n$.


\subsection{$\sigma = 1$ case}


In the $\sigma = 1$ case, $I_n^\sigma (\gvec \alpha, \gvec \gamma)$ has the integral representations
\be
I_n^1 (\gvec \alpha, \gvec \gamma) = \int_{S_K} \! \df^{K - 1} y \, \bigg( \sum_{k = 1}^K \gamma_k y_k \bigg) ^n \prod_{i = 1}^K y_i^{\alpha_i - 1}
\ee
and
\be
I_n^1 (\gvec \alpha, \gvec \gamma) = \fr{1}{\prod_{l = 1}^K \gamma_l^{\alpha_l}} \int_{S_K} \! \df ^{K - 1} x \, \fr{\prod_{i = 1}^K x_i^{\alpha_i - 1}}{(\sum_{j = 1}^K x_j/\gamma_j)^{\alpha_+ + n}} .
\label{In1b}
\ee
The finite sum representation of \eq{Inmultinomial} gives
\be
I_n^1 (\gvec \alpha, \gvec \gamma) = \fr{n! \prod_{l = 1}^K \Gamma(\alpha_l)}{\Gamma(\alpha_+ + n)} \sum_{n_+ = n} \prod_{i = 1}^K \fr{(\alpha_i)_{n_i} \gamma_i^{n_i}}{n_i!} .
\label{In1sum}
\ee


\subsubsection{Carlson $R$-function}


$I_n^1 (\gvec \alpha, \gvec \gamma)$ is, up to multiplication by a multivariate beta function of $\gvec \alpha$, the $R$-function defined by Carlson \cite{carlson} (see \cite{carlsonbook} for a review),
\be
I_n^1 (\gvec \alpha, \gvec \gamma) = \fr{\prod_{i = 1}^K \Gamma(\alpha_i)}{\Gamma(\alpha_+)} R_n (\gvec \alpha, \gvec \gamma) .
\ee
The Carlson $R$-function is defined by the integral representation
\be
R_n (\vec b, \vec z) = \fr{\Gamma (b_+)}{\prod_{k = 1}^K \Gamma (b_k)} \int_{S_K} \! \df^{K - 1} x \, \bigg( \sum_{j = 1}^K z_j x_j \bigg)^n \prod_{i = 1}^K x_i^{b_i - 1} ,
\ee
where $\Real b_i > 0$, and can be analytically continued to more general values of $\vec b$.  Considered as a function of $\vec z$, the Carlson $R$-function is a multivariate generalization of the ${_2}F{_1}$ hypergeometric function.  The Carlson $R$-function can be considered a rewriting of the Lauricella $F_D$-function, but has the advantage that the symmetry of the full permutation group $S_K$ is manifest, whereas standard representations of the $F_D$-function make manifest only symmetry of the permutation group $S_{K - 1}$.  $n$ may be complex-valued, requiring a branch cut, but for our purposes it suffices to consider non-negative integers $n$.  


\subsubsection{Connection to complete homogeneous symmetric polynomials}


Further specializing to the case $\gvec \alpha = \vec 1$, it is well-known (see e.g.\ Theorem 220 of \cite{halipo}) that
\be
I_n^1(\vec 1 , \gvec \gamma) = \fr{n!}{(n + K - 1)!} h_n(\gvec \gamma) ,
\ee
where (with $h_0(\gvec \gamma) = 1$)
\be
h_n(\gvec \gamma) = \sum_{n_+ = n} \prod_{i = 1}^K \gamma_i^{n_i} = \sum_{i_1 \leq i_2 \leq \ldots \leq i_n} \gamma_{i_1} \gamma_{i_2} \ldots \gamma_{i_n}
\ee
is the complete homogeneous symmetric polynomial of degree $n$ in $\gvec \gamma$.  For example, if $K = 3$, then $h_1(\gvec \gamma) = \gamma_1 + \gamma_2 + \gamma_3$ and $h_2(\gvec \gamma) = \gamma_1^2 + \gamma_2^2 + \gamma_3^2 + \gamma_1 \gamma_2 + \gamma_1 \gamma_3 + \gamma_2 \gamma_3$.  The polynomial $h_n(\gvec \gamma)$ is a sum of $\binom{n + K - 1}{n}$ monomials; the corresponding complete homogeneous symmetric mean defined as
\be
q_n(\gvec \gamma) = \binom{n + K - 1}{n}^{-1} h_n(\gvec \gamma) = R_n(\vec 1, \gvec \gamma)
\ee
satisfies $q_n(\vec 1) = 1$, so equivalently
\be
I_n^1(\vec 1 , \gvec \gamma) = \fr{q_n(\gvec \gamma)}{(K - 1)!} .
\ee

We may also express $h_n(\gvec \gamma)$ in terms of power sum symmetric polynomials $p_d(\gvec \gamma) = \sum_{j = 1}^K \gamma_j^d$ as (see e.g.\ equations (7.17) and (7.22) of \cite{stanley})
\be
h_n (\gvec \gamma) = \sum_{\substack{
m_d \geq 0\\
m_1 + 2 m_2 + \ldots + n m_n = n
}} \prod_{d = 1}^n \fr{p_d (\gvec \gamma)^{m_d}}{(m_d)! d^{m_d}} .
\ee
The sum is over all partitions of $n$, so e.g.\ $h_1 = p_1$, $h_2 = \tf{1}{2} (p_1^2 + p_2)$, $h_3 = \tf{1}{6} (p_1^3 + 3 p_1 p_2 + 2 p_3)$.  We may express $h_n$ in terms of $p_1 , \ldots , p_n$ through a complete Bell polynomial as
\be
h_n (p_1 , \ldots , p_n) = \fr{1}{n!} B_n (p_1 , 1! p_2 , 2! p_3 , \ldots , (n - 1)! p_n) .
\label{hbell}
\ee


\subsubsection{Expression in terms of deformed symmetric polynomials}


Returning to general $\gvec \alpha$, we can find an expression for the normalization constant in terms of $\gvec \alpha$-deformed symmetric polynomials,
\be
I_n^1 (\gvec \alpha , \gvec \gamma) = \fr{\prod_{i = 1}^K \Gamma(\alpha_i)}{\Gamma (\alpha_+ + n)} n! h_n (\gvec \alpha, \gvec \gamma) ,
\label{In}
\ee
where we define (with $h_0(\gvec \alpha, \gvec \gamma) = 1$)
\be
h_n (\gvec \alpha, \gvec \gamma) = \sum_{\substack{
m_d \geq 0\\
m_1 + 2 m_2 + \ldots + n m_n = n
}} \prod_{d = 1}^n \fr{p_d (\gvec \alpha, \gvec \gamma)^{m_d}}{(m_d)! d^{m_d}} ,
\label{hn}
\ee
where
\be
p_d (\gvec \alpha, \gvec \gamma) = \sum_{j = 1}^K \alpha_j \gamma_j^d .
\label{pi}
\ee
For practical computation, the sum over partitions in \eq{hn} may be cumbersome.  However, the Newton identities for the standard symmetric polynomials also hold for the $\gvec \alpha$-deformed symmetric polynomials, since $h_n(p_1 , \ldots , p_n)$ continues to be given by \eq{hbell}, and we may still compute $h_n(\gvec \alpha, \gvec \gamma)$ iteratively using
\be
h_n (p_1 , \ldots , p_n) = \fr{1}{n} \sum_{m = 1}^n p_m h_{n - m} (p_1 , \ldots , p_{n - m}) .
\ee

The proof of \eq{In} presented here follows the method suggested by Schl\"{o}milch of differentiating the integral representation \eq{In1b} with respect to $\gvec \gamma^{-1}$.  We first make the induction hypothesis that \eq{In} holds with
\be
h_n (\gvec \alpha, \gvec \gamma) = \sum_{\substack{
m_d \geq 0\\
m_1 + 2 m_2 + \ldots + n m_n = n
}} \fr{\prod_{i = 1}^K p_d (\gvec \alpha, \gvec \gamma)^{m_d}}{z_{n, m_1 \ldots m_n}} ,
\label{hypothesis}
\ee
where the coefficients $z_{n, m_1 \ldots m_n}$ are potentially infinite, but independent of $\gvec \alpha$ and $\gvec \gamma$.  The base $n = 0$ case is the known value of $I_0^1 (\gvec \alpha, \gvec \gamma)$ obtained from the Schl\"{o}milch integral and given in \eq{I0}.  Setting $m = 1$ in the differential recurrence relation \eq{diffrecurm}, we have
\be
I_{n + 1}^1 (\gvec \alpha, \gvec \gamma) = \fr{1}{\alpha_+ + n} \sum_{i = 1}^K \bigg( \gamma_i^2 \fr{\pd}{\pd \gamma_i} + \alpha_i \gamma_i \bigg) I_n^1 (\gvec \alpha , \gvec \gamma) .
\label{diffrecur}
\ee
When differentiating $I_n^1(\gvec \alpha, \gvec \gamma)$, note that
\be
\sum_{j = 1}^K \gamma_j^2 \fr{\pd p_d(\gvec \alpha, \gvec \gamma)}{\pd \gamma_j} = d p_{d + 1} (\gvec \alpha, \gvec \gamma) .
\ee
This implies that $I_{n + 1}(\gvec \alpha, \gvec \gamma)$ is of the form of the induction hypothesis (with $n + 1$ replacing $n$), but with some unknown coefficients, independent of $\gvec \alpha$ and $\gvec \gamma$, multiplying the individual terms appearing in the corresponding expression for $h_{n + 1}(\gvec \alpha, \gvec \gamma)$ in \eq{hn}.

Having established the general form \eq{hypothesis}, we can now fix the coefficients $z_{n, m_1 \ldots m_n}$.  Since these coefficients are independent of $\gvec \alpha$, we can simplify by taking $\gvec \alpha = \vec 1$.  However, we know in this special case that $h_n (\vec 1 , \gvec \gamma)$ and $p_d (\vec 1, \gvec \gamma)$ are complete homogeneous and power sum symmetric polynomials in $\gvec \gamma$ respectively, so $z_{n, m_1 \ldots m_n} = \prod_{d = 1}^n (m_d)! d^{m_d}$, proving \eq{hn} and \eq{pi}.

An alternative proof uses the finite sum representation \eq{In1sum}.  In Lemma 1 of van Laarhoven and Kalker \cite{vanlkal}, it is shown that
\be
\sum_{n_+ = n} \prod_{i = 1}^K \fr{(\alpha_i)_{n_i} \gamma_i^{n_i}}{n_i!} = h_n(\gvec \alpha, \gvec \gamma) ,
\ee
which immediately proves \eq{In}.  The derivation of the lemma relies on manipulations of the generating function \cite{carlson69, menon}
\be
\prod_{i = 1}^K (1 - \gamma_i t)^{-\alpha_i} = \sum_{n = 0}^\infty t^n \fr{(\alpha_+)_n}{n!} R_n(\gvec \alpha, \gvec \gamma) .
\ee
A note in \cite{bullen} states that the generating function is associated with symmetric forms introduced much earlier, however the corresponding references are difficult to access.


\subsection{Applications of simplex integral representation of $h_n$}


We have found that the evaluation of integrals on a simplex results in expressions closely related to complete homogeneous symmetric polynomials.  Although not the focus of our work here, we take the opportunity to remark on the converse viewpoint, that the simplex integrals provide a useful representation of complete homogeneous symmetric polynomials, allowing for natural generalizations.


\subsubsection{Fractional degree complete homogeneous symmetric polynomials}


It has been suggested that the definition of complete homogeneous symmetric polynomials $h_n(X_1 , \ldots , X_K)$ be extended from non-negative integer degree $n$ to an arbitrary real degree $z$ \cite{gaomonyi, bogaomon}.  Three equivalent approaches given are: a 1-dimensional integral based on B-splines, a divided difference expression, and a definition based on Jacobi's bialternant formula.

Complementing these approaches, I define the complete homogeneous symmetric polynomial in a complex-valued degree $z$, for $K \geq 2$, as
\be
h_z (X_1 , \ldots , X_K) = \bigg( \prod_{i = 1}^{K - 1} (z + i) \bigg) \int_{S_K} \! \df^{K - 1} x \, \bigg( \sum_{j = 1}^K X_j x_j \bigg)^z ,
\label{hz}
\ee
which is unambiguous if $z$ is an integer (of any sign) or if all $X_i$ are positive.  More generally the definition depends on a choice of the complex logarithm.  Note that $h_z(X_1, \ldots , X_K) = 0$ for negative integers $z = -1, \ldots , -(K - 1)$.  For $K = 1$, we define $h_z (X_1) = X_1^z$.  Although the representation \eq{hz} involves integrating over more dimensions that the B-spline formula of \cite{gaomonyi, bogaomon}, one advantage is that the integrand is analytic.  The equivalence of the B-spline and simplex integral approaches is a special case of a remarkable relation between multivariate B-splines and Dirichlet averages: see equation (6.3) of \cite{carlson91} specialized to the univariate case, all knots having multiplicity 1 and $g(x) = x^n$.

A closely related integral representation (see e.g.\ \cite{agchgavo} and references within, where it is expressed in terms of expectations involving exponential distributions) is
\be
h_n (X_1 , \ldots , X_K) = \fr{1}{n!} \int_{\bb R_+^K} \! \df^K y \, \expe{-\sum_{i = 1}^K y_i} \bigg( \sum_{j = 1}^K X_j y_j \bigg)^n ,
\label{hz2}
\ee
which reduces to \eq{hz} by setting $y_i = t x_i$ and integrating out $t$.  One advantage of the simplex integral representation \eq{hz} is that the integrand is bounded, whereas the integrand of \eq{hz2} diverges as $\vec y \rightarrow \vec 0$ when $n$ is generalized to negative values.


\subsubsection{Inequalities}


The simplex integral representation \eq{hz} is useful for proving inequalities.  If $n$ is a positive integer, allowing $X_i$ to take any sign, it is manifest that $h_{2 n} (\vec X) \geq 0$, with equality if and only if all $X_i = 0$.  There are a number of different proofs: see \cite{agchgavo} and references within.  The inequality is widely attributed to Hunter \cite{hunter}, who proved a stronger result without using an integral representation.  However, the inequality was proved substantially earlier, appearing in the treatise of Hardy, Littlewood and P\'{o}lya \cite{halipo}, possibly based on work of Schur, using the same integral representation used here.  A more general result of \cite{halipo} is that $\sum_{m = 0}^N \sum_{n = 0}^N q_{m + n}(\vec X) y_m y_n$ is a strictly positive quadratic form for any $N$, since
\be
\sum_{m = 0}^N \sum_{n = 0}^N q_{m + n}(\vec X) y_m y_n = [(K - 1)!]^2 \int_{S_K} \! \df^{K - 1} x \, \bigg[ \sum_{n = 0}^N y_n \bigg( \sum_{i = 1}^K X_i x_i \bigg) ^n \bigg] ^2 \geq 0 .
\label{qinequality}
\ee
Choosing $\vec y$ so that only $y_n$ is non-zero shows that $h_{2 n} (X_1 , \ldots , X_K) \geq 0$.  The general inequality may also be understood by noting that the matrix representing the quadratic form is a Gram matrix with components $Q_{m n} = \langle f_m, f_n \rangle$, where $f_n(\vec X) = (K - 1)! (\sum_{i = 1}^K X_i x_i)^n$ and $\langle f, g \rangle = \int_{S_K} \! \df^{K - 1} x \, f(\vec x) g(\vec x)$.  Since the functions $f_n$ are linearly independent, the corresponding Gram matrix $Q_{m n}$ is positive-definite.  The inequality $q_n (\vec X) ^2 \leq q_{n + m} (\vec X) q_{n - m} (\vec X)$, where $m$ is a positive integer, with equality if and only if all $X_i$ are equal, is also provided in \cite{halipo}, which holds because $| \begin{smallmatrix}
q_{n - m} & q_n \\
q_n & q_{n + m}
\end{smallmatrix} |$ is a minor of a positive-definite matrix.

We may generalize these results to degrees that are not positive integers.  If $n$ is a positive integer and $K < 2 n$ is odd, then the integral representation implies, without any restriction on the sign of $X_i$, that
\be
(-1)^{n - 1} h_{-2 n} (\vec X) \geq 0 .
\ee
Equality holds if and only if all $X_i = 0$.

By the Cauchy--Schwarz inequality, we have, for $X_i \geq 0$, and real $r$ and $s$,
\be
q_{(r + s)/2}(\vec X)^2 \leq q_r (\vec X) q_s(\vec X) ,
\ee
where the complete homogeneous symmetric mean of degree $z$ is
\be
q_z (\vec X) = (K - 1)! \int_{S_K} \! \df ^{K - 1} x \, \bigg( \sum_{j = 1}^K X_j x_j \bigg)^z = \fr{(K - 1)!}{\prod_{i = 1}^{K - 1} (z + i)} h_z(\vec X) . 
\ee
Equality holds if and only if $r = s$ or all $X_i$ are equal.

The Minkowski inequality implies that, for $X_i \geq 0$ and $Y_i \geq 0$, if $p$ is real and $|p| > 1$, then
\be
[h_p (\vec X + \vec Y)]^{1/|p|} \leq [h_p (\vec X)]^{1/|p|} + [h_p (\vec Y)]^{1/|p|} ,
\ee
and if $p$ is real and $|p| < 1$, then
\be
[h_p (\vec X + \vec Y)]^{1/|p|} \geq [h_p (\vec X)]^{1/|p|} + [h_p (\vec Y)]^{1/|p|} .
\ee
In both cases, equality holds if and only if $\vec X$ and $\vec Y$ are linearly dependent or if $\vec X$ or $\vec Y$ vanishes.  This type of inequality has been considered in \cite{agchgavo}, which traces the result to Whiteley \cite{whiteley} and McLeod \cite{mcleod}.

By generalizing the positivity \eq{qinequality} by consideration of $\sum_m \sum_n q_{m + n} (\gvec \alpha, \gvec X)$, Bennett \cite{bennett} has shown that $R_n(\gvec \alpha, \vec z)$, considered as a sequence indexed by non-negative integers $n$ and assuming that $\alpha_i$ are positive, is a Hamburger moment sequence for arbitrary real $\vec z$, i.e.\ $R_n(\gvec \alpha, \vec z) = \int_{-\infty}^\infty \! \df \mu(x) \, x^n$ for some positive measure $\df \mu$.  Moreover,  if $z_i$ are positive, then it is a Stieltjes moment sequence, i.e.\ $R_n(\gvec \alpha, \vec z) = \int_0^\infty \! \df \mu(x) \, x^n$ for some positive measure $\df \mu$.  For the case $\gvec \alpha = \vec 1$, an explicit measure is provided by the univariate B-spline with $K$ simple knots located at $z_i$,
\be
F(x; \vec z) = (K - 1) \sum_{i = 1}^K \fr{(z_i - x)_+^{K - 2}}{\prod_{j \neq i} (z_i - z_j)} = \fr{K - 1}{2} \sum_{i = 1}^K \fr{| z_i - x | (z_i - x)^{K - 3}}{\prod_{j \neq i} (z_i - z_j)} ,
\ee
where $x_+ = \max (0, x) = \tf{1}{2} (x + |x|)$, since it is well-known that the B-spline $F(x; \vec z)$ is non-negative and has support on $(\min (\vec z) , \max (\vec z))$.  This gives \cite{neuman}
\be
q_n (\vec z) = \int_{-\infty}^\infty \! \df x \, x^n F(x; \vec z) .
\ee
One may similarly show that $I_n^\sigma(\gvec \alpha, \gvec \gamma)$ is a Hamburger moment sequence, and moreover a Stieltjes moment sequence if all $\gamma_i > 0$.


\section{Properties of the distributions}



\subsection{Exponential family}


If all parameters are allowed to vary, then the Dirichlet mixture distribution $DM(\gvec \alpha, \gvec \gamma, n, \sigma)$ and, more generally, the Schl\"{o}milch mixture distribution $SM(\gvec \alpha, \gvec \beta, \gvec \gamma, n, \sigma, \tau)$ do not belong to the exponential family of distributions.  However, if all parameters except $\gvec \alpha$ are fixed, then they are exponential family distributions, i.e.\ the probability density functions can be expressed in the form
\be
f(\vec x) = c(\gvec \alpha) h(\vec x) \exp \bigg( \sum_{i = 1}^K \alpha_i T_i (\vec x) \bigg) .
\ee
The natural parameter is $\gvec \alpha$, and $\vec T(\vec X)$ is a sufficient statistic.  General moments for exponential family distributions are given by
\be
\bb E (T_1^{\ell_1} \ldots T_K^{\ell_K}) = c(\gvec \alpha) \fr{\pd^{\ell_+}}{\pd \alpha_1^{\ell_1} \ldots \pd \alpha_K^{\ell_K}} \fr{1}{c(\gvec \alpha)} ,
\ee
where $\ell_+ = \sum_{i = 1}^K \ell_i = || \gvec \ell ||_1$.  In particular, the means are
\be
\bb E (T_i) = - \fr{\pd}{\pd \alpha_i} \log c(\gvec \alpha) ,
\ee
and the covariances are
\be
\Cov (T_i, T_j) = - \fr{\pd^2}{\pd \alpha_i \, \pd \alpha_j} \log c(\gvec \alpha) .
\ee

For the Schl\"{o}milch mixture distribution, with $f(\vec x)$ given by \eq{SM}, the sufficient statistic is
\be
\vec T (\vec X) = \tau (\log X_1 , \ldots , \log X_K) - \log \bigg( \sum_{j = 1}^K \beta_j X_j^\tau \bigg) \vec 1 ,
\ee
and we may take
\begin{align}
h(\vec x) & = \fr{\tau^{K - 1} (\sum_{k = 1}^K \beta_k^\sigma \gamma_k x_k^{\sigma \tau})^n}{(\sum_{j = 1}^K \beta_j x_j^\tau)^{\sigma n} \prod_{i = 1}^K x_i} , & c(\gvec \alpha) & = \fr{\prod_{l = 1}^K \beta_l^{\alpha_l}}{I_n^\sigma(\gvec \alpha, \gvec \gamma)} .
\end{align}


\subsection{Log-ratios}


The sum constraint implies that $2 \sum_{j \neq i} \Cov (X_i, X_j) = -\Var (X_i)$, so there is a tendency for $\Cov(X_i, X_j)$ to be negative.  For this reason, use of log-ratios $\log (X_i/X_j)$ have been suggested as more appropriate for analysis of compositional data \cite{aitchison}, with the log-ratio covariance $\Cov [\log(Y_i/Y_j), \log(Y_k/Y_l)]$ used to examine correlations.


\subsubsection{Dual random variables on the simplex}


Suppose that $\vec X$ and $\vec Y$ are random variables on the simplex $S_K$, related by the duality of \eq{duality}, i.e.
\begin{align}
Y_i & = \fr{\beta_i X_i^\tau}{\sum_{j = 1}^K \beta_j X_k^\tau} , & X_i & = \fr{(Y_i/\beta_i)^{1/\tau}}{\sum_{j = 1}^K (Y_j/\beta_j)^{1/\tau}} ,
\end{align}
where $\beta_i > 0$ and $\tau > 0$, repeated here for convenience.  The duality between $\vec X$ and $\vec Y$ leads to simple relations between their log-ratio means and covariances.

There are relations of expectations given by
\be
\bb E \bigg[ \log \bigg( \fr{X_i^\tau}{\sum_{k = 1}^K \beta_k X_k^\tau} \bigg) \bigg] = \bb E \bigg[ \log \bigg( \fr{Y_i}{\beta_i} \bigg) \bigg] 
\ee
and
\begin{align}
& \bb E \bigg[ \log \bigg( \fr{X_i^\tau}{\sum_{k = 1}^K \beta_k X_k^\tau} \bigg) \log \bigg( \fr{X_j^\tau}{\sum_{k = 1}^K \beta_k X_k^\tau} \bigg) \bigg] = \bb E \bigg[ \log \bigg( \fr{Y_i}{\beta_i} \bigg) \log \bigg( \fr{Y_j}{\beta_j} \bigg) \bigg] .
\end{align}
It follows that the log-ratio means are related by
\be
\bb E \bigg[ \log \bigg( \fr{X_i}{X_j} \bigg) \bigg] = \fr{1}{\tau} \bb E \bigg[ \log \bigg( \fr{Y_i}{Y_j} \bigg) \bigg] - \fr{1}{\tau} \log \bigg( \fr{\beta_i}{\beta_j} \bigg) ,
\ee
and the log-ratio covariances are related by
\begin{align}
& \Cov \bigg[ \log \bigg( \fr{X_i}{X_j} \bigg) , \log \bigg( \fr{X_k}{X_l} \bigg) \bigg] = \fr{1}{\tau^2} \Cov \bigg[ \log \bigg( \fr{Y_i}{Y_j} \bigg) , \log \bigg( \fr{Y_k}{Y_l} \bigg) \bigg] .
\end{align}


\subsubsection{Dirichlet and Schl\"{o}milch mixture distributions}


As a particular example of the duality, if $\vec X \sim SM(\gvec \alpha, \gvec \beta, \gvec \gamma, n, \sigma, \tau)$ and $\vec Y \sim DM(\gvec \alpha, \gvec \gamma, n, \sigma)$, then we obtain the expectations
\be
\bb E \bigg[ \log \bigg( \fr{\beta_i X_i^\tau}{\sum_{k = 1}^K \beta_k X_k^\tau} \bigg) \bigg] = \bb E (\log Y_i) = \fr{\pd}{\pd \alpha_i} \log I_n^\sigma(\gvec \alpha, \gvec \gamma)
\ee
and
\begin{align}
& \bb E \bigg[ \log \bigg( \fr{\beta_i X_i^\tau}{\sum_{k = 1}^K \beta_k X_k^\tau} \bigg) \log \bigg( \fr{\beta_j X_j^\tau}{\sum_{k = 1}^K \beta_k X_k^\tau} \bigg) \bigg] = \bb E ( \log Y_i \, \log Y_j ) \nnr
& = \fr{\pd^2}{\pd \alpha_i \, \pd \alpha_j} \log I_n^\sigma(\gvec \alpha, \gvec \gamma) + \bigg( \fr{\pd}{\pd \alpha_i} \log I_n^\sigma(\gvec \alpha, \gvec \gamma) \bigg) \bigg( \fr{\pd}{\pd \alpha_j} \log I_n^\sigma(\gvec \alpha, \gvec \gamma) \bigg) .
\end{align}
It follows that the log-ratio means are given by
\be
\tau \bb E \bigg[ \log \bigg( \fr{X_i}{X_j} \bigg) \bigg] + \log \bigg( \fr{\beta_i}{\beta_j} \bigg) = \bb E \bigg[ \log \bigg( \fr{Y_i}{Y_j} \bigg) \bigg] = \bigg( \fr{\pd}{\pd \alpha_i} - \fr{\pd}{\pd \alpha_j} \bigg) \log I_n^\sigma(\gvec \alpha, \gvec \gamma) ,
\label{ElogY}
\ee
and the log-ratio covariances are given by
\begin{align}
& \tau^2 \Cov \bigg[ \log \bigg( \fr{X_i}{X_j} \bigg) , \log \bigg( \fr{X_k}{X_l} \bigg) \bigg] = \Cov \bigg[ \log \bigg( \fr{Y_i}{Y_j} \bigg) , \log \bigg( \fr{Y_k}{Y_l} \bigg) \bigg] \nnr
& = \bigg( \fr{\pd}{\pd \alpha_i} - \fr{\pd}{\pd \alpha_j} \bigg) \bigg( \fr{\pd}{\pd \alpha_k} - \fr{\pd}{\pd \alpha_l} \bigg) \log I_n^\sigma (\gvec \alpha, \gvec \gamma) .
\end{align}
Using derivatives computed in Section \ref{lown}, we can obtain explicit expressions for low $n$.

For $n = 0$, from \eq{d20} we obtain
\be
\tau^2 \Cov \bigg[ \log \bigg( \fr{X_i}{X_j} \bigg) , \log \bigg( \fr{X_k}{X_l} \bigg) \bigg] = \psi'(\alpha_i) (\delta_{i k} - \delta_{i l}) - \psi'(\alpha_j) (\delta_{j k} - \delta_{j l}) .
\ee
In particular, if $i, j, k, l$ are distinct, then $\Cov [ \log (X_i/X_j), \log (X_k/X_l)] = 0$, so these log-ratios are uncorrelated.  If $j = l$ and $i, j, k$ are distinct, then
\be
\tau^2 \Cov \bigg[ \log \bigg( \fr{X_i}{X_j} \bigg) , \log \bigg( \fr{X_k}{X_j} \bigg) \bigg] = \psi'(\alpha_j) . 
\ee
For positive $x$, the trigamma function $\psi'(x)$ is positive, so $\log (X_i/X_j)$ and $\log (X_k/X_j)$ must have positive correlation.

For $n = 1$, from \eq{d21} we obtain
\begin{align}
\tau^2 \Cov \bigg[ \log \bigg( \fr{X_i}{X_j} \bigg) , \log \bigg( \fr{X_k}{X_l} \bigg) \bigg] & = \bigg( p_i \psi'(\alpha_i + \sigma) + (1 - p_i) \psi'(\alpha_i) + \fr{c_i^2}{p_i} \bigg) (\delta_{i k} - \delta_{i l}) \nnr
& \qquad - \bigg( p_j \psi'(\alpha_j + \sigma) + (1 - p_j) \psi'(\alpha_j) + \fr{c_j^2}{p_j} \bigg) (\delta_{j k} - \delta_{j l}) \nnr
& \qquad - (c_i - c_j) (c_k - c_l) .
\end{align}
If $i, j, k, l$ are distinct, then $\Cov [ \log (X_i/X_j), \log (X_k/X_l)] = - (c_i - c_j) (c_k - c_l)/\tau^2$ is non-zero, unlike the $n = 0$ case.  If $j = l$ and $i, j, k$ are distinct, then
\be
\tau^2 \Cov \bigg[ \log \bigg( \fr{X_i}{X_j} \bigg) , \log \bigg( \fr{X_k}{X_j} \bigg) \bigg] = p_j \psi'(\alpha_j + \sigma) + (1 - p_j) \psi'(\alpha_j) + \fr{c_j^2}{p_j} - (c_i - c_j) (c_k - c_j) .
\ee
Specializing to a Dirichlet mixture, for which $\tau = 1$, this recovers the log-ratio covariance given in \cite{ongmig}.

For $n = 1$, further specializing to the case that $\sigma = 1$, we have
\be
\tau^2 \Cov \bigg[ \log \bigg( \fr{X_i}{X_j} \bigg) , \log \bigg( \fr{X_k}{X_l} \bigg) \bigg] = \psi'(\alpha_i) (\delta_{i k} - \delta_{i l}) - \psi'(\alpha_j) (\delta_{j k} - \delta_{j l}) - \fr{(\gamma_i - \gamma_j) (\gamma_k - \gamma_l)}{h_1 (\gvec \alpha, \gvec \gamma)^2} ,
\ee
so $\Cov [ \log (X_i/X_j), \log (X_k/X_l)] = - (\gamma_i - \gamma_j) (\gamma_k - \gamma_l)/\tau^2 h_1 (\gvec \alpha, \gvec \gamma)^2$ if $i, j, k, l$ are distinct.  If $j = l$ and $i, j, k$ are distinct, then
\be
\tau^2 \Cov \bigg[ \log \bigg( \fr{X_i}{X_j} \bigg) , \log \bigg( \fr{X_k}{X_j} \bigg) \bigg] = \psi'(\alpha_j) - \fr{(\gamma_i - \gamma_j) (\gamma_k - \gamma_j)}{h_1 (\gvec \alpha, \gvec \gamma)^2} ,
\ee
which one can easily see may be positive or negative.

For $n \geq 2$, there is similarly non-zero correlation between $\log(X_i/X_j)$ and $\log(X_k/X_l)$, and positive or negative correlation between $\log(X_i/X_j)$ and $\log(X_k/X_j)$.


\subsection{Moments}
\label{moments}


Mixed moments $\bb E (\prod_{i = 1}^K X_i^{\ell_i})$ of the Dirichlet distribution, and therefore of mixtures of Dirichlet distributions, can be explicitly computed straightforwardly.  In general, moments for the Schl\"{o}milch distribution $S(\gvec \alpha, \gvec \beta, \tau)$ or the more general Schl\"{o}milch mixture distribution $SM(\gvec \alpha, \gvec \beta, \gvec \gamma, n, \sigma, \tau)$ involve integrals that cannot be expressed in terms of standard analytic functions.


\subsubsection{Explicitly tractable cases}


For mixed moments $\bb E (\prod_{i = 1}^K X_i^{\ell_i})$ with $\ell_i$ given by non-negative integers, a special case that can be dealt with more explicitly is $SM(\gvec \alpha, \gvec \beta, \overline{\gvec \beta^{-1}}, n, 1, 1)$, which includes (by taking $n = 0$) the scaled Dirichlet or G3D distribution $S(\gvec \alpha, \gvec \beta, 1)$, and the tilted Dirichlet distribution $SM(\gvec \alpha, \gvec \alpha, \overline{\gvec \alpha^{-1}}, 1, 1, 1)$.  The probability density function, from \eq{case3} with $\tau = 1$, is
\be
f(\vec x) = \fr{( \prod_{l = 1}^K \beta_l^{\alpha_l} )  \prod_{i = 1}^K x_i^{\alpha_i - 1}}{I_n^1 (\gvec \alpha, \gvec \beta^{-1}) (\sum_{j = 1}^K \beta_j x_j)^{\alpha_+ + n}} ,
\label{tilt}
\ee
which we henceforth assume in this section.

Using the duality relation for $I_n^\sigma(\gvec \alpha, \gvec \gamma)$ of \eq{Idual}, the moments can be expressed as 
\be
\bb E ({\textstyle \prod_{i = 1}^K X_i^{\ell_i}}) = \fr{I_{n - \ell_+}^1 (\gvec \alpha + \gvec \ell, \gvec \beta^{-1})}{(\prod_{i = 1}^K \beta_i^{\ell_i}) I_n^1 (\gvec \alpha, \gvec \beta^{-1})} ,
\ee
or equivalently as
\be
\bb E (X_{k_1} \ldots X_{k_{\ell_+}}) = \fr{I_{n - \ell_+}^1 (\gvec \alpha + \sum_{a = 1}^{\ell_+} \vec e_{k_a}, \gvec \beta^{-1})}{(\prod_{a = 1}^{\ell_+} \beta_{k_a}) I_n^1 (\gvec \alpha, \gvec \beta^{-1})} .
\ee
Analogous expressions for moments of the $S$-distribution of Dickey \cite{dickey83} are expressed in terms of the Carlson $R$-function, but without use of the function's duality relation.

Recall that, if $n$ is a non-negative integer, the function $I_n^1 (\gvec \alpha, \gvec \beta^{-1})$ can be expressed using \eq{In} in terms of explicit $\gvec \alpha$-deformed symmetric polynomials.  Therefore, the lowest-order moments, for $\ell_+ \leq n$, are given by closed-form expressions without any need to evaluate integrals.  For $n \geq 1$, the first moments are
\be
\bb E (X_i) = \fr{\alpha_i h_{n - 1}(\gvec \alpha + \vec e_i, \gvec \beta^{-1})}{n \beta_i h_n(\gvec \alpha, \gvec \beta^{-1})} .
\ee
For $n \geq 2$, the second moments are
\be
\bb E (X_i X_j) = \fr{(\alpha_i \alpha_j + \alpha_i \delta_{i j}) h_{n - 2}(\gvec \alpha + \vec e_i + \vec e_j, \gvec \beta^{-1})}{n (n - 1) \beta_i \beta_j h_n(\gvec \alpha, \gvec \beta^{-1})} ,
\ee
and the covariances, for $i \neq j$, are
\begin{align}
\Cov (X_i, X_j) & = \fr{\alpha_i \alpha_j}{\beta_i \beta_j n^2 (n - 1) h_n(\gvec \alpha, \gvec \beta^{-1})^2} [n h_{n - 2}(\gvec \alpha + \vec e_i + \vec e_j, \gvec \beta^{-1}) h_n(\gvec \alpha, \gvec \beta^{-1}) \nnr
& \qquad - (n - 1) h_{n - 1}(\gvec \alpha + \vec e_i, \gvec \beta^{-1}) h_{n - 1}(\gvec \alpha + \vec e_j, \gvec \beta^{-1})].
\end{align}

Expressions for moments with $\ell_+ > n$ are more complicated, but $I_{n - \ell_+}^1 (\gvec \alpha + \gvec \ell, \gvec \beta^{-1})$, and hence the moments, can be reduced to one-dimensional integrals.  By using the Schl\"{o}milch integral \eq{schlomilchint} on $S_2$ and $S_K$, we compute
\begin{align}
& \int_{S_K} \! \df^{K - 1} x \, \fr{\prod_{i = 1}^K x_i^{\alpha_i + \ell_i - 1}}{(\sum_{j = 1}^K \beta_j x_j)^{\alpha_+ + n}} \nnr
& = \fr{\Gamma (\alpha_+ + \ell_+)}{\Gamma(\alpha_+ + n) \Gamma(\ell_+ - n)} \int_{S_K} \! \df^{K - 1} x \, \int_0^1 \! \df u \, \fr{u^{\alpha_+ + n- 1} (1 - u)^{\ell_+ - n- 1} \prod_{i = 1}^K x_i^{\alpha_i + \ell_i - 1}}{[(\sum_{j = 1}^K \beta_j x_j) u + 1 - u]^{\alpha_+ + \ell_+}} \nnr
& = \fr{\Gamma (\alpha_+ + \ell_+)}{\Gamma(\alpha_+ + n) \Gamma(\ell_+ - n)} \int_0^1 \! \df u \, \int_{S_K} \! \df^{K - 1} x \, \fr{u^{\alpha_+ + n- 1} (1 - u)^{\ell_+ - n- 1} \prod_{i = 1}^K x_i^{\alpha_i + \ell_i - 1}}{[\sum_{j = 1}^K (\beta_j u + 1 - u) x_j]^{\alpha_+ + \ell_+}} \nnr
& = \fr{\prod_{i = 1}^K \Gamma (\alpha_i + \ell_i)}{\Gamma(\alpha_+ + n) \Gamma(\ell_+ - n)} \int_0^1 \! \df u \, \fr{u^{\alpha_+ + n - 1} (1 - u)^{\ell_+ - n - 1}}{\prod_{j = 1}^K (\beta_j u + 1 - u)^{\alpha_j + \ell_j}} ,
\end{align}
assuming that $\ell_i > 0$ and $\ell_+ > n$.  Since rescaling allows us to assume that $\beta_+ = 1$, implying that $0 < \beta_i < 1$, the singularities of the integrand at $u = 1/(1 - \beta_i)$ do not lie on $[0, 1]$.  We obtain the moments, for $\ell_+ > n$,
\be
\bb E ({\textstyle \prod_{i = 1}^K X_i^{\ell_i}}) = \fr{\prod_{i = 1}^K \Gamma (\alpha_i + \ell_i)}{\Gamma(\alpha_+ + n) \Gamma(\ell_+ - n) I_n^1 (\gvec \alpha, \gvec \beta^{-1})} \int_0^1 \! \df u \, \fr{u^{\alpha_+ + n - 1} (1 - u)^{\ell_+ - n - 1}}{\prod_{j = 1}^K (\beta_j u + 1 - u)^{\alpha_j + \ell_j}} .
\ee
For the scaled Dirichlet or G3D distribution $S(\gvec \alpha, \gvec \beta, 1)$ of \eq{schlomilchdist}, given by $n = 0$, these moments have been given in \cite{dickey}.  Using integration by parts, we obtain an analytic continuation that is valid for $\ell_+ > 0$,
\be
\bb E ({\textstyle \prod_{i = 1}^K X_i^{\ell_i}}) = \fr{\prod_{i = 1}^K \Gamma (\alpha_i + \ell_i)}{\Gamma(\alpha_+ + n) \Gamma(\ell_+) I_n^1 (\gvec \alpha, \gvec \beta^{-1})} \int_0^1 \! \df u \, (1 - u)^{\ell_+ - 1} \fr{\df^n}{\df u^n} \bigg( \fr{u^{\alpha_+ + n - 1}}{\prod_{j = 1}^K (\beta_j u + 1 - u)^{\alpha_j + \ell_j}} \bigg) .
\ee


\subsubsection{$n = 1$ case}


For $n = 1$, using $I_1^1(\gvec \alpha, \gvec \beta^{-1})$ given by \eq{I11}, the probability density function of \eq{tilt} is
\be
f(\vec x) = \fr{(\prod_{l = 1}^K \beta_l^{\alpha_l}) \Gamma (1 + \alpha_+)}{(\prod_{j = 1}^K \Gamma(\alpha_j)) \sum_{k = 1}^K \alpha_k/\beta_k} \fr{\prod_{i = 1}^K x_i^{\alpha_i - 1}}{(\sum_{j = 1}^K \beta_j x_j)^{\alpha_+ + 1}} .
\ee
The first moments are
\be
\bb E (X_i) = \fr{\alpha_i}{\beta_i h_1 (\gvec \alpha, \gvec \beta^{-1})} = \fr{\alpha_i/\beta_i}{\sum_{j = 1}^K \alpha_j/\beta_j} .
\ee

Although $\Cov(X_i, X_j)$ is always negative for the Dirichlet distribution, numerical computation indicates that it can generally have either sign in this $n = 1$ case.  Moreover, the correlation $\Corr(X_1, X_2)$ can be arbitrarily close to 1, as shown by numerical computation with $K = 3$, $\gvec \alpha = (1/\epsilon, 1/\epsilon, 1)$, $\gvec \beta = (1, 1, \epsilon)/(2 + \epsilon)$, taking the limit $\epsilon \rightarrow 0$. 


\subsubsection{$n = 2$ case}


We may similarly study the $n = 2$ case of the probability density function of \eq{tilt}.  A benefit of this case is the tractability of its second moments, which can be computed in closed-form.  The first moments are
\be
\bb E (X_i) = \fr{\alpha_i h_1(\gvec \alpha + \vec e_i, \gvec \beta^{-1})}{2 \beta_i h_2(\gvec \alpha, \gvec \beta^{-1})} ,
\ee
where
\begin{align}
h_1 (\gvec \alpha + \vec e_i, \gvec \beta^{-1}) & = p_1 (\gvec \alpha + \vec e_i, \gvec \beta^{-1}) = \fr{1}{\beta_i} + \sum_{k = 1}^K \fr{\alpha_k}{\beta_k} , \nnr
2 h_2 (\gvec \alpha, \gvec \beta^{-1}) & = p_1(\gvec \alpha, \gvec \beta^{-1})^2 + p_2(\gvec \alpha, \gvec \beta^{-1}) = \bigg( \sum_{k = 1}^K \fr{\alpha_k}{\beta_k} \bigg)^2 + \sum_{k = 1}^K \fr{\alpha_k}{\beta_k^2} .
\end{align}
The second moments are
\be
\bb E (X_i X_j) = \fr{(\alpha_i \alpha_j + \alpha_i \delta_{i j})}{2 \beta_i \beta_j h_2(\gvec \alpha, \gvec \beta^{-1})} ,
\ee
which gives the variance
\be
\Var(X_i) = \fr{\alpha_i [2 (\alpha_i + 1) h_2(\gvec \alpha, \gvec \beta^{-1}) - \alpha_i h_1(\gvec \alpha + \vec e_i, \gvec \beta^{-1})^2]}{4 \beta_i^2 h_2 (\gvec \alpha, \gvec \beta^{-1})^2}
\ee
and, for $i \neq j$, the covariance
\be
\Cov (X_i , X_j) = \fr{\alpha_i \alpha_j [2 h_2(\gvec \alpha, \gvec \beta^{-1}) - h_1(\gvec \alpha + \vec e_i, \gvec \beta^{-1}) h_1(\gvec \alpha + \vec e_j, \gvec \beta^{-1})]}{4 \beta_i \beta_j h_2(\gvec \alpha, \gvec \beta^{-1})^2} .
\ee

Although $\Cov(X_i, X_j)$ is always negative for the Dirichlet distribution, it can generally have either sign for this $n = 2$ case, determined by the sign of
\be
2 h_2 (\gvec \alpha, \gvec \beta^{-1}) - h_1 (\gvec \alpha + \vec e_i, \gvec \beta^{-1}) h_1 (\gvec \alpha + \vec e_j, \gvec \beta^{-1}) = \sum_{k = 1}^K \fr{\alpha_k}{\beta_k^2} - \fr{1}{\beta_i \beta_j} - \bigg( \fr{1}{\beta_i} + \fr{1}{\beta_j} \bigg) \sum_{k = 1}^K \fr{\alpha_k}{\beta_k} .
\ee
For example, if $K = 3$, $\alpha_i = 1$, $\beta_3/\beta_1 \ll 1$ and $\beta_3/\beta_2 \ll 1$, then $\Cov(X_1, X_2)$ will be positive.  Moreover, the correlation $\Corr(X_1, X_2)$ can be arbitrarily close to 1, as shown by the limit $\epsilon \rightarrow 0$ in the example with $\gvec \alpha = (1/\epsilon, 1/\epsilon, 1)$ and $\gvec \beta = (1, 1, \epsilon)/(2 + \epsilon)$, for which $\Corr(X_1, X_2) = (1 - 4 \epsilon - \epsilon^2)/(1 + 6 \epsilon + \epsilon^2)$. 


\subsection{Distributions on superellipsoids}


A simple generalization is to transform the simplex coordinates by $x_i \rightarrow (1 - c_i) (x_i/b_i)^{a_i}$ for some $a_i > 0$, $b_i > 0$ and $c_i < 1$.  We could set $c_i = 0$ without loss of generality, but it is useful to include for comparison with literature and for taking a $c_i \rightarrow 1$ limit.  The simplex then becomes the positive orthant of a $(K - 1)$-dimensional generalized superellipsoid, $E_{K, \vec a, \vec b, \vec c}^+ = \{ \vec x = (x_1 , \ldots , x_K)  \subset \bb R^K : \sum_{i = 1}^K (1 - c_i) (x_i/b_i)^{a_i} = 1, x_i \geq 0 \}$.

For simplicity and for comparison with literature, we start with the Schl\"{o}milch distribution $S(\gvec \alpha, \gvec \beta, 1)$, given by \eq{schlomilchdist}, and set $\beta_i = \beta_K/(1 - c_i)$ for $i = 1 , \ldots , K - 1$.  Transforming the simplex to a generalized superellipsoid, we find that
\be
f(\vec x) = \fr{\Gamma(\alpha_+)}{\prod_{l = 1}^K \Gamma(\alpha_l)} \bigg( \prod_{k = 1}^{K - 1} \fr{a_k x_k^{a_k \alpha_k - 1}}{b_k^{a_k \alpha_k}} \bigg) \fr{[1 - \sum_{i = 1}^{K - 1} (1 - c_i) (x_i/b_i)^{a_i}]^{\alpha_K - 1}}{[1 + \sum_{j = 1}^{K - 1} c_j (x_j/b_j)^{a_j}]^{\alpha_+}}
\ee
is a probability distribution function on $E_{K, \vec a, \vec b, \vec c}^+$, expressed in terms of independent coordinates $x_1 , \ldots , x_{K - 1}$.  This matches the multivariate generalized beta (MGB) distribution of Cockriel and McDonald \cite{cocmcd}, who also allow for negative $a_i$, although they do not explicitly check the normalization of the probability distribution function.

The MGB distribution is significant because it includes a large number of well-known distributions as special cases.  Its specialization to a univariate distribution is the 5-parameter generalized beta (GB) distribution of McDonald and Xu \cite{mcdxu}.  The authors have found applications to economics and finance, but special cases are well-known and have applications in many fields. These include distributions with support on a finite interval, such as the beta and Pareto distributions, and with support on a half-line, such as the generalized gamma, inverted beta, and Lomax distributions.  The distributions presented in this article lead to further unification of distributions in the literature.


\section{Discussion}


Whereas the Dirichlet distribution on the $(K - 1)$-dimensional simplex $S_K$ has $K$ independent continuous parameters, the more general Schl\"{o}milch mixture distribution given by the probability density function \eq{SM} has $3K$ independent continuous parameters and a discrete parameter $n$, providing considerably more freedom.  This unifies a number of families of distributions on the simplex that have been studied in the literature.  Unlike the Dirichlet distribution, the Schl\"{o}milch mixture distribution allows for positive covariances $\Cov(X_i, X_j)$, non-zero log-ratio covariance $\Cov[\log (X_i/X_j), \log(X_k/X_l)]$, negative $\Cov[\log(X_i/X_j), \log(X_k/X_j)]$, and multi-modality.  We have also constructed an inverse Schl\"{o}milch distribution that generalizes the Concrete or Gumbel-softmax distribution.

One could consider distributions in which the non-negative integer parameter $n$ appearing in the Schl\"{o}milch mixture distribution is generalized to be continuous.  There would, however, be additional complications, because of a lack of identifiability of $n$ in the limit of a Dirichlet distribution, and because the method for generating samples for integer $n$ from gamma distributions does not generalize to non-integer $n$.

For the subfamily of Schl\"{o}milch mixture distributions with $\sigma = 1$, we have found that the normalization constant is closely related to complete homogeneous symmetric polynomials.  Conversely, our study motivates a definition, through an $S_K$ simplex integral representation, of complete homogeneous symmetric polynomials in which the degree is not restricted to a non-negative integer.  Such integral representations may be useful in deriving further inequalities for symmetric functions.

\end{document}